\magnification=\magstep1
\vbadness=10000
\hbadness=10000
\tolerance=10000

\def\C{{\bf C}}           
\def\Q{{\bf Q}}           
\def\R{{\bf R}}           
\def\Z{{\bf Z}}           
\def\e{{\bf e}}           
\def\Tr{{\rm Tr}}         
\def\sign{{\rm sign}}     
\def\Aut{{\rm Aut}}       

\proclaim Reflection groups of Lorentzian lattices. 
\hfill 21 September 1999.

Richard E. Borcherds, 
\footnote{$^*$}{ 
Supported by a Royal Society professorship and NSF grant DMS-9970611. 
This paper was partly written at the Max-Planck institute in Bonn.}

Math Dept, 
Evans Hall \#3840,
University of California at Berkeley, 
CA 94720-3840, 
U.S.A. 

e-mail: reb@math.berkeley.edu
home page: http://math.berkeley.edu/\hbox{\~{}}reb

\bigskip

\proclaim Contents.

1. Introduction.

Notation.

2.~Modular forms. 

3.~Discriminant forms and the Weil representation. 

4.~The singular theta correspondence. 

5.~Theta functions. 

6.~Eta quotients. 

7.~Dimensions of spaces of modular forms. 

8.~The geometry of $\Gamma_0(N)$. 

9.~An application of Serre duality.

10.~Eisenstein series. 

11.~Reflective forms. 

12. Examples. 

13.~Open problems. 

\proclaim
1.~Introduction.

The aim of this paper is to provide evidence for the following
new principle: interesting reflection groups of Lorentzian lattices are
controlled by certain modular forms with poles at cusps. We use this
principle to explain many of the known examples of such reflection
groups, and to find several new examples of reflection
groups of Lorentzian lattices, including one whose fundamental domain
has 960 faces.

We do not give a precise definition of what it means for a reflection
group of a Lorentzian lattice to be interesting, mainly because there
seem to be occasional counterexamples to almost any precise version of
the principle.  However the interesting groups should include the
cases when the reflection group is cofinite, or more generally when
the quotient of the full automorphism group by the reflection group
contains a free abelian subgroup of finite index.

The main idea of this paper is roughly as follows (and is described 
in more detail in section 11). Suppose that
$L$ is a Lorentzian lattice in $\R^{1,n}$. Then the idea is that if
$L$ has an interesting reflection group, there should be a modular form
of weight ${1-n\over 2}$ and level $N$ with certain rather mild singularities
(called reflective singularities) at cusps. The singular theta correspondence
of [B98a]
associates a piecewise linear function on the hyperbolic space of $\R^{1,n}$
to this modular form, and the singularities of this piecewise linear function 
should be reflection hyperplanes of the reflection group. 

We do not have any sort of proof that interesting Lorentzian
reflection groups always correspond to reflective modular forms (or
even a good definition of what an interesting group is). However we
can still use this rather vague correspondence to find new reflection
groups, because it is usually easy to list examples of reflective
forms using the results in the first 10 sections. In section 12 we
use this to
give many examples of reflection groups corresponding to reflective
modular forms. In particular we show that most of the known examples
in dimensions at least 5 can be found by systematically searching for
reflective forms of small levels. (This does not seem to work so well
for Lorentzian lattices of dimension 3 and possibly 4; Nikulin
has found many examples that do not obviously correspond to reflective
forms.) This gives some sort of structure to the collection 
of all nice Lorentzian reflection groups, which previously looked 
like a miscellaneous collection of unrelated examples found by 
half a dozen assorted methods. We also find many new examples 
of Lorentzian lattices with interesting reflection groups, including
one whose fundamental domain has 960 faces (comfortably 
beating the previous record
of 210 faces). The main problem with these examples is that there 
are almost too many of them: for small composite levels (4, 6, 8, 9) there
are so many cases that we do not try to list them all but just give
vague indications of how to list many of them. Somewhere round about 
level 25 the number of examples for each level decreases to a trickle, 
and above this level we only seem to get a few low dimensional lattices
for each level. However there are probably still a few occasional examples 
when the level is a few hundred. 

So reflective modular forms seem to be a useful practical method for
finding interesting reflection groups of Lorentzian lattices. On the
other hand there are occasional counterexamples to show that it does
not always work. There are one or two high dimensional reflection
groups which seem well behaved but do not appear to correspond to
reflective modular forms, and on the other hand there are occasional
Lorentzian lattices with non-zero reflective modular forms whose
reflection groups are quite complicated. Moreover, when there is a
reflective form, the reflection group of the Lorentzian lattice can
have several different behaviors: for examples it might be cofinite,
or it might have a free abelian subgroup of finite index, or it might
have a space-like Weyl vector, or it might have none of these
properties. So the relation between reflection group and reflective
forms is useful in practice, but the theoretical side is still rather
mysterious. The reader who wishes to tidy up the theory is warned that 
section 12 contains counterexamples to several plausible 
simplifying conjectures.

Sections 2 to 10 are mainly a summary of various assorted results that
we need, most of which are minor variations of known results. The
contents of these sections should mostly be clear from their titles.

Correction. 
A. Barnard pointed out a mistake in the formula for the Weyl vector 
in [B98a, theorem 10.4]. In the formula for $\rho_z$, the terms for 
$\lambda =0$ were incorrectly omitted. So the condition 
$(\lambda,W)>0$ in the sum should be deleted, and the factor
of ${1\over 2}$ in front of the sum should be replaced by ${1\over 4}$. 

I would like to thank 
D. Allcock, 
E. Freitag, 
S. Kondo, 
I. Dolgachev, 
and D. Zagier 
for their help. 

\proclaim Notation and terminology.

\item{$\tilde{}$} 
A metaplectic double cover of a group. 
\item{$\sqrt{}$} 
The principal value of the square root, with $-\pi/2<\arg(\sqrt{*})\le \pi/2$.
\item{${}'$}
The dual of a lattice. 
\item{$A$} 
A discriminant form. 
\item{$A_n$} 
The $A_n$ root lattice, or the elements of order $n$ in a discriminant
form $A$.
\item{$A^n$} 
The $n$'th powers of elements of $A$.
\item{$\Gamma$} 
A discrete subgroup of $Mp_2(\R)$. 
\item{$\Delta$}
The Dedekind delta function $\eta^{24}$. 
\item{$D_n$} 
The $D_n$ root lattice. 
\item{$e_\gamma$} 
An element of a basis of $\C[L'/L]$.
\item{$\e$} 
$\e(x)=\exp({2\pi i x})$, $\e_n(x)=\exp({2\pi i x/n})$. 
\item{$E_k$} 
An Eisenstein series (see section 10),
or the $E_k$ root lattice. 
\item{$\eta$}
The Dedekind eta function. 
\item{$\theta_L$} 
A theta function of a lattice $L$. 
\item{$g$} 
The genus of a subgroup of $SL_2(\R)$. 
\item {$II_{m,n}$} 
The even unimodular lattice of dimension
$m+n$ and signature $m-n$.
\item{$L$} 
An even lattice. 
\item{$N$} 
The level of a modular form or discriminant form. 
\item {$q^n$} 
$e^{2\pi i n\tau}$ or a discriminant form. 
\item {$\Q$} 
The rational numbers.
\item {$\R$} 
The real numbers.
\item {$\rho_L$} 
A representation of $Mp_2(\Z)$.
\item{$R_{j}$ }
The primitive elliptic element fixing the point $j$. 
\item{$\sign$}
The signature of a lattice or discriminant form. 
\item{$SL$} 
A special linear group.
\item {$\tau$} 
A complex number with positive imaginary part, or the number of
orbits of cusps. 
\item{$Tr$} 
The trace of something. 
\item{$T_{a/c}$}
The primitive parabolic element fixing the cusp $a/c$. 
\item {$\Z$} 
The integers.
\item{$\chi$} 
A character. For $\chi_\theta$, $\chi_n$, see section 5. 

\proclaim 2.~Modular forms.

In this section we recall the definition of a vector valued  modular form
and set up notation for the rest of the paper.

We define $\e(x)$ to be $\exp(2\pi i x)$, and we define $\e_n(x)$ to be
$\exp(2\pi x/n)$

Recall that the group $SL_2(\R)$ has a (metaplectic) 
double cover $\widetilde{SL}_2(\R)$,
whose
elements can be written in the form
$$\left(\pmatrix{a&b\cr c&d\cr},\pm\sqrt{c\tau +d}\right)
$$
where $\pmatrix{a&b\cr c&d\cr}\in SL_2(\R)$.
The multiplication is defined so that the usual formulas for
the transformation of modular forms work for half integer weights,
which means that
$$
(A,f(\cdot))(B,g(\cdot))= (AB,f(B(\cdot))g(\cdot))
$$
for $A,B\in SL_2(\R)$ and $f,g$ suitable functions on $H$.
The group $\widetilde{SL}_2(\Z)$ is the inverse image of
$SL_2(\Z)\subset SL_2(\R)$ in $\widetilde{SL}_2(\R)$.
The group $\widetilde{SL}_2(\Z)$ is generated by $S$ and $T$
where $T=\left(\pmatrix{1&1\cr 0&1\cr},1\right)$ and $S=\left(\pmatrix{0&-1\cr
1&0\cr},\sqrt{\tau}\right)$,
with $S^2=(ST)^3=Z$, $Z=\left(\pmatrix{-1&0\cr 0&1\cr},i\right)$.
The center is cyclic of order and 4 is generated by $Z$. 
The quotient by $Z^2$ is the group 
$SL_2(\Z)$.

Suppose that $\Gamma$ is a discrete subgroup of $\widetilde{SL}_2(\R)$ 
that contains $Z$ and is co-finite (this means that the quotient 
space has finite volume). 
Suppose that $\rho$ is a representation of $ \Gamma$
on a finite dimensional complex vector space $V_\rho$.
Choose $k\in \Q$.
We define a modular form of weight $k$ and type
$\rho$ 
to be a holomorphic function $f$ on the upper half plane $H$
with values in the vector space $V_\rho$ such that
$$f\left({a\tau+b\over c\tau+d}\right)
= (c\tau+d)^k
\rho\left(\left(\pmatrix{a&b\cr c&d\cr},\sqrt{c\tau+d}\right)\right)
f(\tau)
$$
for elements $\left(\pmatrix{a&b\cr c&d\cr},\sqrt{c\tau+d}\right)$
of $\Gamma$. The expression $(c\tau+d)^k$
means of course $\sqrt{c\tau+d}^{2k}$ with the principal value of
$\sqrt{}$. 
(We allow singularities at cusps.)

A modular form has a Fourier expansion at the cusp at infinity as follows.
The Fourier coefficients $c_{n,\gamma}\in \C$ of $f$ are defined by
$$f(\tau) = \sum_{n\in \Q}\sum_{\gamma} c_{n,\gamma}q^ne_\gamma$$
where $q^n$ means $\e( n\tau)$ and where the sum runs over a
basis $e_\gamma$ of $V_\rho$ consisting of eigenvectors of $T$. Note
that $n$ is not necessarily integral; more precisely, $c_{n,\gamma}$
is nonzero only if $n\equiv \lambda_\gamma\bmod 1$, where the
eigenvalue of $T$ on $e_\gamma$ is $\e( \lambda_\gamma)$. We
say that $f$ is meromorphic at the cusp $i\infty$ if $c_{n,\gamma}=0$ for
$n<<0$, and we say $f$ is meromorphic at the cusp $a/c$ if
$f((a\tau+b)/(c\tau+d))$ is meromorphic at $i\infty$ for 
$\pmatrix{a&b\cr c&d\cr}\in SL_2(\Z)$. We say that $f$ is holomorphic 
at cusps if the coefficients of the Fourier expansions
at all cusps vanish for $n<0$. 

\proclaim 3.~Discriminant forms and the Weil representation.

In this section we recall the definition of the Weil representation 
of a discriminant form, and prove some results about it that will be
used in section 11. 

We let $L$ be a nonsingular even lattice of dimension $\dim(L)$ and signature
$\sign(L)$, with dual
$L'$. The quotient $L'/L$ is a finite group whose order
is the absolute value of the discriminant of the lattice $L$. 
We define a discriminant form $A$ to be a finite abelian group 
with a $\Q/\Z$-valued quadratic form $\gamma\mapsto \gamma^2/2$. 
If $L$ is an even lattice then $L'/L$ is a discriminant form, 
with the quadratic form of $L'/L$ given by the mod 1 reduction of
$\gamma^2/2$, 
and conversely every discriminant form can be constructed in this way. 
(For the theory of the discriminant form of
a lattice see [N].)
This quadratic form on $L'/L$ determines the signature mod 8 of $L$, 
by Milgram's formula
$$
\sum_{\gamma\in L'/L} \e(\gamma^2/2)
= \sqrt{|L'/L|}\e(\sign(L)/8).
$$
We define the signature $\sign(A)\in \Z/8\Z$ 
of a discriminant form to be the signature
mod 8
of any even lattice with that discriminant form. 
We let the elements $e_\gamma$ for $\gamma\in L'/L$ be the standard
basis of the group ring $\C[L'/L]$, so that $e_\gamma
e_\delta=e_{\gamma+\delta}$. 

A particularly important example $\rho_A$ of a 
unitary representation of $\widetilde{SL}_2(\Z)$,
called the Weil representation of the discriminant form $A$, 
can be constructed as follows. 
The underlying space of $\rho_A$ is the group ring
$\C[A]$ of $A$, and the action is defined by
$$\rho_A(T)(e_\gamma) = \e{((\gamma,\gamma)/2)}e_\gamma
$$
$$\rho_A(S)(e_\gamma) = {\e(-\sign(A)/8)\over \sqrt{|A|}}
\sum_{\delta\in A} \e(-(\gamma,\delta))e_\delta
$$
where $S$ and $T$ are the standard generators of $\widetilde{SL}_2(\Z)$. 
The representation $\rho_A$
factors through the double cover $\widetilde{SL}_2(\Z/N\Z)$ of the 
finite group $SL_2(\Z/N\Z)$, where $N$ is a positive integer such that
$N\gamma^2/2$ is an integer for all
$\gamma\in L'$. The smallest such integer $N$ is called the
level of $A$. In particular the representation $\rho_A$
factors through a finite quotient of $\widetilde{SL}_2(\Z)$.
If $L$ is an even lattice then we define $\rho_L$ to be the
representation $\rho_{L'/L}$.

We summarize some results about discriminant forms from [C-S, chapter
15, section 7]; for more details see [N] or [C-S]. We use a minor
variation of the notation of [C-S] for discriminant forms. We recall
that every discriminant form can be written as a sum of Jordan
components (not uniquely if $p=2$), and every Jordan component can be
written as the sum of indecomposable Jordan components (usually not
uniquely). The possible non-trivial Jordan components are as follows.
We let $q>1$ be a power of a prime $p$ and $n$ a positive integer and
$t\in \Z/8\Z$. We define $antisquare$ by $antisquare(q^{\pm n})=0$ if
$q$ is a square or the exponent is $+n$, and $antisquare(q^{\pm n})=1$
if $q$ is not a square and the exponent is $-n$. (See [C-S] page 370.)

For $q$ odd the non-trivial Jordan components 
of exponent $q$ are $q^{\pm n}$ for
$n\ge 1$.  The indecomposable
components are $q^{\pm 1}$, generated by an element $\gamma$
with $q\gamma=0$, $\gamma^2\equiv a/q\bmod 2$ where $a$ is an even
integer with ${a\choose p}=\pm 1$. The component $q^{\pm n}$ is a sum
of copies of $q^{+1}$ and $q^{-1}$, with an even number of copies of
$q^{-1}$ if $\pm n=+n$ and an odd number if $\pm n=-n$.
These components all have level $q$.
The signature is given by $\sign(q^{\pm n})=-n(q-1)+4antisquare(q^{\pm n})$.

For $q$ even the odd Jordan components of exponent $q$
are $q^{\pm n}_t$. If $n=1$
then $t\equiv \pm 1\bmod 8$ if $\pm=+$ and $t\equiv \pm 3\bmod 8$ if
$\pm=-$. If $n=2$ then $t\equiv 0,\pm 2\bmod 8$ if $\pm=+$ and
$t\equiv 4,\pm 2\bmod 8$ if $\pm=-$. For any $n$ we have $t\equiv
n\bmod 2$. The indecomposable components are $q^{\pm 1}_t$ for
${t\choose 2}=\pm 1$ and are generated by an element $\gamma$ with
$q\gamma=0$, $\gamma^2\equiv t/q\bmod 2$. (Note that some of these
are isomorphic to each other.) These components all have level $2q$. 
The signature is given by $\sign(q^{\pm n}_t)=t+4antisquare(q^{\pm n})$.

For $q$ even the non-trivial
even Jordan components of exponent $q$ are $q^{\pm 2n}=q^{\pm 2n}_{II}$. 
The indecomposable even Jordan components are $q^{\pm 2}$, 
which are generated by 2 elements $\gamma$ and $\delta$ 
with $q\gamma=q\delta=0$, $(\gamma,\delta)=1/q$, 
$\gamma^2\equiv\delta^2\equiv0\bmod 2$ 
if $\pm =+$, $\gamma^2\equiv\delta^2\equiv 2/q \bmod 2$ if $\pm =-$. 
These components all have level $q$. 
The signature is given by $\sign(q^{\pm n})=4antisquare(q^{\pm n})$.

The sum of two Jordan components with the same prime power $q$ can be
worked out as follows: we add the ranks, multiply the signs in the
exponent, and if any components have a
subscript $t$ we add together all subscripts $t$.

If $A$ is a discriminant form, 
then we define $A_n$ to be the elements of order $n$. 
We define $A^n$ to be the $n$'th powers of elements of $A$, 
so that we have an exact sequence
$$0\longrightarrow A_n\longrightarrow A \longrightarrow A^n
\longrightarrow 0$$ and $A^n$ is the orthogonal complement of $A_n$.
We define $A^{n*}$ to be the set of elements $\delta\in A$ such that
$(\gamma,\delta)\equiv n\gamma^2/2\bmod 1$ for all $\gamma\in A_n$, so
that $A^{n*}$ is a coset of $A^n$. We easily see that $A^n$ is the
same as $A^{n*}$ if and only if the Jordan block of type $2^k$ (where $2^k||n$)
is even. In any case, $A^{n*}$ always contains an element $\delta$ with
$2\delta=0$. 

\proclaim Lemma 3.1. Suppose that $A$ is a discriminant form. 
Then 
$$\sum_{\gamma\in A} \e((\gamma,\delta)-n\gamma^2/2)$$
is 0 unless $\delta\in A^{n*}$ (in which case it has absolute value
$\sqrt{|A||A_n|}$). 

Proof. 
The square of the absolute value of this sum is
$$
\eqalign{
&\sum_{\gamma_1,\gamma_2\in A} 
\e((\gamma_1,\delta)-n\gamma_1^2/2-(\gamma_2,\delta)+n\gamma_2^2/2)\cr
=&
\sum_{\gamma_1,\gamma_2\in A} 
\e((\gamma_1,\delta)-n\gamma_1^2/2-n(\gamma_1,\gamma_2))\cr
=&
|A|\sum_{\gamma_1\in A_n} 
\e((\gamma_1,\delta)-n\gamma_1^2/2)\cr
}$$
The map taking $\gamma_1$ to $ \e((\gamma_1,\delta)-n\gamma_1^2/2)$
is a character of $A_n$, so this sum is 0 unless this is the trivial 
character, in other words unless $\delta\in A^{n*}$. 
This proves lemma 3.1. 

\proclaim Lemma 3.2. 
Suppose that $g\in \widetilde{SL}_2(\Z)$ has image $\pmatrix{a&b\cr c&d\cr}\in 
SL_2(\Z)$. Then $\rho_A(e_0)$ is a linear combination of the elements
$e_\gamma$ for $\gamma\in A^{c*}$. 

Proof. It is sufficient to prove this when $g$ is of the form 
$T^mST^nS$ for some $m,n\in \Z$ with $(N,n)=(N,c)$, 
because $g$ is a product of an element of this form 
by an element of $\tilde\Gamma_0(N)$ (where $N$ is the level of $A$)
and $e_0$ is an eigenvalue of $\tilde\Gamma_0(N)$.

We calculate the image of $1=e_0\in \C[A]$ under these elements of
$\widetilde{SL}_2(\Z)$. We know that
$$
S(e_0) 
= {\e(-\sign(A)/8)\over \sqrt{|A|}}
\sum_{\gamma\in L'/L} e_\gamma
$$
Applying $T^n$ shows that
$$
T^{n}S(1)
= {\e(-\sign(A)/8)\over \sqrt{|A|}}
\sum_{\gamma\in A}\e(n(\gamma,\gamma)/2) e_\gamma.
$$
Applying $S$ again shows that
$$
ST^{n}S(1)
= {\e(-\sign(A)/4)\over {|A|}}
\sum_{\delta\in A}\sum_{\gamma\in A}
\e((\gamma,\delta)+n(\gamma,\gamma)/2) e_\delta
$$

Using lemma 3.1, we see that
the coefficient of $e_\delta$ in this expression 
is 0 unless $\delta\in A^{c*}$. 
As all the elements $e_\delta$ are eigenvectors of $T^m$, the same is true
for $T^mST^nS(e_0)$.
This proves lemma 3.2.

\proclaim 4.~The singular theta correspondence.

We summarize some of the results from [B98a]. The main idea is that we
can use modular forms with poles at cusps to construct some
automorphic forms with singularities. In particular we can often use
this to construct piecewise linear functions on hyperbolic space with
singularities along the reflection hyperplanes of a reflection group,
and this gives the connection between modular forms with singularities 
and nice hyperbolic reflection groups. 

If $L$ is a lattice then we define the Grassmannian $G(L)$
to be the set of maximal positive definite subspaces of $L\otimes \R$. 
It is a symmetric space acted on by the orthogonal group $O_L(\R)$. 

The Siegel theta function
$\theta_{L+\gamma} $ of a coset $L+\gamma$ of $L$ in $L'$ is defined by
$$\theta_L(\tau;v^+)= \sum_{\lambda\in L+\gamma}\e(\tau\lambda_{v^+}^2/2
+\bar\tau\lambda_{v^-}^2/2)
$$
for $\tau\in H$, $v^+\in G(L)$.
We will write
$\Theta_L$ for the $\C[L'/L]$-valued function
$$\Theta_L(\tau;v)=\sum_{\gamma\in
L'/L}e_\gamma\theta_{L+\gamma}(\tau;v)
{}.
$$
Siegel's transformation formula for $\Theta_L$ under $\widetilde{SL}_2(\Z)$
([B98a, theorem 4.1])
is given by
$$
\Theta_L\left({a\tau+b\over c\tau+d}v\right)
= (c\tau+d)^{b^+/2}(c\bar\tau+d)^{b^-/2}
\rho_L\left(\pmatrix{a&b\cr c&d\cr},\sqrt{c\tau+d} \right)
\Theta_L(\tau;v).
$$

We define $\Phi(v,F)$ by 
$$
\Phi(v,F)=\int_{SL_2\backslash H}\bar\Theta_L(\tau; v)F(\tau)y^{b^+/2-2}dxdy
$$
as in section 6 of [B98a]. By theorem 6.2 of [B98a], $\Phi(v,F)$ is
an automorphic function of $v\in G(L)$ whose only singularities 
are on points of the form $\gamma^\perp$, for $\gamma\in L'$, 
$\gamma^2<0$, where there
is a nonzero coefficient $c_{\gamma^2/2,\gamma}$ of $F$. 

\proclaim Theorem 4.1. Suppose $L$ is an even lattice of signature
$(2,b^-)$ and $F$ is a modular form of weight $1-b^-/2$ and
representation $\rho_L$ which is holomorphic on $H$ and meromorphic at cusps
and whose
coefficients $c_\lambda(m)$ are integers for $m\le 0$.
Then there is a meromorphic function $\Psi_L(Z_L,F)$
for $Z\in P$ with the following properties.
\item{1.} $\Psi_L(Z_L,F)$ is an automorphic form of weight $c_0(0)/2$
for the group $\Aut(L,F)$ with respect to some unitary character $\chi$
of $\Aut(L,F)$.
\item{2.} The only zeros or poles of $\Psi_L$ lie on the
rational quadratic divisors $\lambda^\perp$
for $\lambda\in L$, $\lambda^2<0$ and are zeros of order
$$\sum_{0<x\in \R\atop x\lambda\in L'}c_{x\lambda}(x^2\lambda^2/2).
$$
(or poles if this number is negative).
\item{3.} $\Psi_L$ is a holomorphic function
if the orders of all zeros in item 2 above are nonnegative.
If in addition $L$ has
dimension at least 5, or if $L$ has dimension 4 and contains no 2 dimensional
isotropic sublattice, then $\Psi_L$ is a holomorphic automorphic form.
If in addition $c_0(0)=b^--2$ then $\Psi_L$ has singular weight 
so the only nonzero Fourier
coefficients of $\Psi_L$ correspond to vectors of $K$ of norm 0.

This follows from theorem 13.3 of [B98a]. 

If $L$ is Lorentzian, in other words if $\sign(L)=2-\dim(L)$, 
then the set of all 1 dimensional positive definite subspaces of
$L$ is a copy of hyperbolic space of dimension $\dim(L)-1$. 

\proclaim Theorem 4.2. Suppose $M$ is a Lorentzian lattice of dimension
$1+b^-$. Suppose that $F$ is a modular form of type $\rho_M$
and weight $(1/2-b^-/2,0)$ which is holomorphic on $H$ and meromorphic
at cusps and all of whose Fourier coefficients $c_\lambda(m)$
are real for $m<0$.
Finally suppose that if $c_\lambda(\lambda^2/2)\ne 0$ and
$\lambda^2<0$ then reflection in $\lambda^\perp$ is in $\Aut(M,F,C)$. Then
$\Aut(M,F,C)$ is the semidirect product of a reflection subgroup
and a subgroup fixing the Weyl vector $\rho(M,W,F)$ of a Weyl chamber $W$.

This is a special case of theorem 12.1 of [B98a].

Both theorems 4.1 and 4.2 depend on integrating the vector valued
modular form against a vector valued theta function over a fundamental
domain of $SL_2(\Z)$.  In this paper we usually start with a complex
valued modular form for $\Gamma_0(N)$ rather than a vector valued form
as used in theorems 4.1 and 4.2. There are two more or less equivalent
ways to use these theorems on complex valued forms of level
$N$. First, instead of integrating a vector valued form times the
vector valued theta function of a lattice over a fundamental domain of
$SL_2(\Z)$, we can integrate a scalar valued modular form times the
theta function of a lattice over a fundamental domain of
$\Gamma_0(N)$. Alternatively we can first induce the complex valued
modular form for $\Gamma_0(N)$ up to a vector valued modular form for
$SL_2(\Z)$, and then apply the theorems directly to part of this
vector valued form. For these constructions to work, it is necessary
and sufficient for the complex valued form to be a modular form for
some character $\chi$ of $\tilde\Gamma_0(N)$, where the scalar valued
theta function of the lattice is a modular form of character $\chi$
and level $\sign(L)/2$ for $\tilde\Gamma_0(N)$.  Several sections of
this paper describe how to find such modular forms. Note that the
singularities of the automorphic form associated to a level $N$
modular form depend on all poles at all cusps of this form, not just
the poles at $i\infty$.

Theorem 4.2 is very useful in practice for finding Lorentzian
lattices with interesting reflection groups, because we just find
lattices together with modular forms satisfying the conditions of the theorem.
However there is a problem with using it for theoretical purposes: it
seems hard to give useful general conditions under which the Weyl
vector is nonzero or has positive norm. If the Weyl vector happens to
be zero then of course theorem 4.2 does not say anything. In practical
examples this does not matter because we can just check in each case
to see whether the vector is zero (which does happen occasionally).
Note the rather curious fact that in this paper we do not need to use
the fact that theorem 4.2 has been proved (or even that it is true!)
because we are only using it to suggest interesting places to look for
lattices, and whenever we find a lattice using theorem 4.2 we still have 
to prove its properties directly because of the possibility that
the Weyl vector is 0. 

\proclaim 5.~Theta functions.

In this section we work out the level and character of theta functions
of even lattices. Most of the results are known, but there seems to
be no convenient reference giving the results in the generality we
require.

\proclaim Lemma 5.1. 
Suppose that $N$ is a positive integer. 
If $4\not|N$ then two characters of $\Gamma_0(N)$ are the same provided that 
they have the same values 
on the elements 
such that $c>0$, $d>0$, and $d\equiv 1\bmod 4$. 
If $4|N$ then two characters of $\tilde\Gamma_0(N)$ are the same provided that 
they have the same values on $Z$ 
and on the elements 
such that $c>0$, $d>0$, and $d\equiv 1\bmod 4$. 

Proof. 
It is sufficient to show that the images of 
the elements mentioned above generate $\Gamma_0(N)$. 
Suppose that $\pmatrix{a&b\cr c&d\cr}\in \Gamma_0(N)$. We will show how to 
multiply it by powers of elements of the generating set above so
that it becomes an element of the generating set, which will prove the lemma. 
We first note that $T$ is in
the group generated by the set above, because the generating set is
closed under left multiplication by $T$ and is nonempty. If $d$ is
even then $c$ is odd so we can multiply it on the right by $T$ so that
$d$ is odd, hence we can assume that $d$ is odd.

Next we arrange that $d\equiv 1\bmod 4$. If $4|N$ we multiply 
by $Z$ if necessary so that $d\equiv 1\bmod 4$. 
If $4\not|N$ we multiply on the right by $\pmatrix{1&0\cr N&1\cr}$
if necessary to make $c$ not divisible by $4$, and then multiply 
on the right by a suitable power of $T$ so that $d\equiv 1\bmod 4$. 

We now have to make $c$ and $d$ positive (without changing $d\bmod 4$). 
We multiply 
on the right by a suitable power of $\pmatrix{1&0\cr 4N&1\cr}$ to make $c$
positive without changing $d$. 
Finally we multiply 
on the right by a suitable power of $\pmatrix{1&4\cr 0&1\cr}$ to make $d$
positive without changing $c$.
The result is in the generating set, so this proves lemma 5.1.

We define the symbol ${c\choose d}$ for all pairs of coprime 
integers $c$ and $d$ as follows. The symbol is multiplicative in both
$c$ and $d$. If $d$ is an odd prime it is just the usual Legendre
symbol. If $d$ is 2 it is 1 if $c\equiv \pm 1\bmod 8$ and $-1$
otherwise. Finally if $d=-1$ it is 1 if $c>0$ and $-1$ if $c<0$. 
Finally we define ${0\choose \pm 1}={\pm 1\choose 0}=1$.

We now define some characters $\chi_n$ (for $n$ a positive integer)
and $\chi_\theta$ of $\tilde \Gamma_0(N)$. 
We suppose that if $p$ is an odd prime occurring an
odd number of times in the prime factorization of $n$ then it divides
$N$. Also suppose that if $2$ occurs an odd number of times in the
prime factorization of $n$ then $8$ divides $N$. We define the
character $\chi_n$ of $\Gamma_0(N)$ by
$$\chi_n\left(\pmatrix{ a&b\cr c&d\cr}\right)
= {d\choose n}.
$$

\proclaim Lemma 5.2. 
Let $\theta_{A_1}=\sum_{n\in Z}q^{n^2}$ be the theta function of
the $A_1$ lattice. 
There is a (unique) character $\chi_\theta$ of
the metaplectic double cover of $\Gamma_0(4)$ such 
that 
$$ \theta_{A_1}\left({a\tau+b\over c\tau+d}\right)
= \chi_\theta\left(\pmatrix{a&b\cr c&d\cr},\sqrt{c\tau+d}\right) 
\sqrt{c\tau+d} \theta_{A_1}(\tau).$$
(In other words $\theta_{A_1}$ is a modular form for $\tilde\Gamma_0(4)$
of weight $1/2$ and character $\chi_\theta$.)
The values are given by
$$
\chi_\theta\left(\pmatrix{a&b\cr c&d\cr},\pm\sqrt{c\tau+d}\right) 
= 
\cases{
\pm {c\choose d}
& if $d\equiv 1\bmod 4$\cr
\pm(-i) {c\choose d}
& if $d\equiv 3\bmod 4$\cr
}
$$
In particular, $\chi_\theta(Z) = -i$.

Proof. This follows from the theorem on page 148 of [Ko]. 

\proclaim Lemma 5.3. Suppose $4|N$. 
Then the kernel of the character $\chi_\theta$ of $\widetilde{\Gamma}_0(N)$
maps isomorphically onto $\Gamma_1(4)\cap \Gamma_0(N)$, 
and if we identify the kernel
with this image then $\widetilde{\Gamma}_0(N)$ is the product
$(\Gamma_1(4)\cap\Gamma_0(N))\times \Z/4\Z$ 
(where $\Z/4\Z$ is its center, generated by
$Z$). The lifting of $ \Gamma_1(4)\cap \Gamma_0(N)$ to
$\tilde \Gamma_0(N)$ is given by 
$$
\pmatrix{a&b\cr c&d\cr}\mapsto
\left(\pmatrix{a&b\cr c&d\cr},{c\choose d}\sqrt{c\tau+d}\right)
$$

Proof. This follows immediately from lemma 5.2 because $\chi_\theta$
is a character whose values are $\pm 1$ or $\pm i$, and 
$\chi_\theta(Z)=-i$, and $\Gamma_0(N) $ is the product of its center of order 2
(generated by $Z$)
and the subgroup
$\Gamma_1(4)\cap \Gamma_0(N)$.
This proves lemma 5.3. 

We define the group $\Gamma_0^2(N)$ to be the subgroup of
$\Gamma_0(N)$ of elements whose diagonal entries are squares in
$\Z/n\Z$. If \hbox{$4\not|N$} then $\Gamma_0^2(N)$ is the intersection of the
kernels of the characters $\chi_p$ for $p$ an odd prime dividing
$N$. If $4|N$ then $\Gamma_0^2(N)$ can be lifted to a subgroup of
$\tilde \Gamma_0(N)$ as in lemma 5.3, and is the intersection of the
kernels of the characters $\chi_\theta$ and $\chi_p$ of $\tilde
\Gamma_0(N)$ for $p$ a prime dividing $N/4$.

\proclaim Theorem 5.4. Suppose that $A$ is a discriminant form of level 
dividing $N$. If $b$ and $c$ are divisible by $N$ then
$g=\left(\pmatrix{a&b\cr c&d\cr}, \sqrt{c\tau+d}\right)
\in \widetilde{SL}_2(\Z)$ acts on
the Weil representation $\C[A]$ by
$$g(e_\gamma) = \chi_A(g)e_{a\gamma}$$ 
where 
is the character
of $\tilde \Gamma_0(N)$ given by
$$
\chi_A = 
\cases{
\chi_\theta^{\sign(A)+{-1\choose |A|}-1}\chi_{|A|2^{\sign(A)}}
& if $4|N$\cr
\chi_{|A|}
& if $4\not| N$\cr
}
$$

Proof. 
First assume that $A$ has even signature. 
Choose an even lattice in a positive
definite space with discriminant form $A$. Then 
[E94, corollary 3.1] and the discussion on [E94, page 94]
show
that 
$$\pmatrix{a&b\cr c&d\cr}(e_\gamma) 
={(-1)^{\sign(A)/2}|A|\choose d} e_{a\gamma}$$
provided that $d$ is odd and positive. 
By lemma 5.1 it is sufficient to check that this is equal to the value
of $\chi_A$ when $g=Z$ and when $d\equiv 1\bmod 4$, $d>0$.
But in the latter case ${-1\choose d}=1$ and ${|A|\choose d} = {d\choose |A|}$,
so this has the same character values as $\chi_{|A|}$. 
Also if $d\equiv 1\bmod 4$ and $\sign(A)$ is even then
$\chi_{2^{\sign(A)}}$ 
and $\chi_\theta^{\sign(A)}$ and $\chi_\theta^{{-1\choose A}-1}$
are all 1 on the element $g$. Therefore the two characters coincide
on elements with $d\equiv 1\bmod 4$. 

As $Z(e_\gamma)=(-i)^{\sign(A)}e_{-\gamma}$ we see that
$\chi_A(Z)=(-i)^{\sign(A)}$
We now check that the characters are equal on the element $Z$.  If
$4\not|N$ this follows from
$\chi_{|A|}(Z)={-1\choose |A|}=(-1)^{\sign(A)/2}$. If $4|N$
this follows from $\chi_2(Z)=1$,
$\chi_\theta^{\sign(A)}(Z)= (-i)^{\sign(A)}$, and
$\chi_{|A|}(Z)= {-1\choose |A|} = (-i)^{1-{-1\choose |A|}}=
\chi_\theta^{1-{-1\choose |A|}}(Z)$. This proves theorem
5.4 when $A$ has even signature. 

We can do the case of odd signature very quickly by reducing it to the
case of even signature as follows. 
If $A$ has odd signature then $4|N$, and 
the discriminant form $A\oplus\langle 2\rangle$
has and determinant $2|A|$
(where $\langle 2\rangle$ is the discriminant form of
the $A_1$ lattice and has order 2).
Theorem 5.4 for the element $\gamma\in A$ now follows from 
lemma 5.2 and theorem 5.4 applied to the element
$\gamma+0\in A\oplus\langle 2\rangle$. This proves theorem 5.4. 

\proclaim 6.~Eta quotients.

In this section we work out the levels and characters of some
eta quotients. We will use these results in section 12 to construct
examples of modular forms of given characters.

The function $\eta(t\tau)$ has a zero of order $(t,c)^2/24t$ at the cusp $a/c$.

\proclaim Lemma 6.1. (Rademacher.) 
Recall that $\eta(\tau)=q^{1/24}\prod_{n>0}(1-q^n)$ 
is the Dedekind eta function.
Suppose that $\pmatrix{a&b\cr c&d\cr}\in SL_2(\Z)$ with
$c>0$. Then
$$\eta\left({a\tau+b\over c\tau+d}\right)
=
\chi_\eta\left(\pmatrix{a&b\cr c&d\cr}, \sqrt{c\tau+d}\right)
\sqrt{c\tau+d}\eta(\tau)
$$
where $\chi_\eta$ is a character of $\widetilde{SL}_2(\Z)$ 
with values given as follows:
$$
\eqalign{
&\chi_\eta\left(\!\pmatrix{a&\!b\cr c&\!d\cr}, \pm\sqrt{c\tau+d}\right)
=
\cases{
\pm{d\choose c} \e_{24}(-3c+bd(1-c^2)+c(a+d))
&$c$ odd, $c>0$\cr
\pm{-d\choose -c} \e_{24}(3c-6+bd(1-c^2)+c(a+d))
&$c$ odd, $c<0$\cr
\pm{c\choose d} \e_{24}(3d-3+ac(1-d^2)+d(b-c))
&$d$ odd, $c\ge 0$\cr
\pm{-c\choose d} \e_{24}(-3d-9+ac(1-d^2)+d(b-c))\!
&$d$ odd, $c<0$\cr
}
\cr}
$$

Proof. The cases with $c\ge0$ follow from the theorem in 
[R, p. 163]. The cases with $c<0$
follow easily from the cases with $c>0$
and the fact that $\chi_\eta(Z)=-i$.
This proves lemma 6.1. 

\proclaim Theorem 6.2. Suppose that we are given 
a positive integer $N$,
and integers $r_\delta$ for $\delta|N$ and $|A|$
with $|A|/\prod_{\delta|N} \delta^{r_\delta}$ a rational square. 
Suppose that 
${1\over 24}\sum_{\delta|N} r_\delta\delta$ and
${N\over 24}\sum_{\delta|N} r_\delta/\delta$ are both integers.
Then 
$$\prod_{\delta|N}\eta(\delta\tau)^{r_\delta}$$ is a modular form for
$\tilde\Gamma_0(N)$ of weight $k=\sum_\delta r_\delta/2$ and character
equal to $\chi_{|A|}$ if $4\not|N$, and to
$\chi_\theta^{2k+{-1\choose |A|}-1}\chi_{ 2^{2k}|A|}$ if $4|N$. 

Proof. By lemma 5.1 it is enough to check that the characters
are equal whenever $c>0$ and $d\equiv 1\bmod 4$, and also that
they are equal on $Z$ if $4|N$.

If $c>0$ and $d\equiv 1\bmod 4$ then by lemma 6.1 the 
character value is given by
$$
\eqalign{
&\prod_{\delta|N} \left(\e_{24}(3d-3+a(c/\delta)(1-d^2)+d(b\delta-c/\delta))
{c/\delta\choose d}\right) ^{r_\delta}
\cr
=&
\e_{24}\left(db\sum_{\delta|N} \delta r_\delta
+(a-d-ad^2)c\sum_{\delta|N} r_\delta/\delta
+3(d-1)\sum_{\delta|N} r_\delta \right)
{c\choose d}^{\sum_{\delta|N}{r_\delta}}
{\prod_{\delta|N}\delta^{r_\delta} \choose d}
\cr
=& i^{ (d-1)k}
{c\choose d}^{2k}
{|A| \choose d}
\cr
=&{d\choose 2}^{2k}
{c\choose d}^{2k}
{d\choose |A| }
\cr
}$$
If $4\not|N$ then $2k$ is even so this is the value of the character
$\chi_{|A|}$. If $4|N$ then this is the value of
$\chi_\theta^{2k}\chi_{2^{2k}|A|}$, which is the same as the value of
$\chi_\theta^{2k+{-1\choose |A|}-1}\chi_{2^{2k}|A|}$ because
$\chi_\theta=\pm 1$ whenever $d\equiv 1\bmod 4$. 
So both characters have the same value whenever $c>0$ and $d\equiv 1\bmod 4$. 

Finally we have to check that both characters are equal on $Z$ whenever 
$4|N$. This follows by the same argument used in theorem 5.4. 
This proves theorem 6.2. 

Theorem 6.2 generalizes some theorem of Newman ([N57], [N59]), who did the
case of weight 0 and trivial character.

\proclaim 7.~Dimensions of spaces of modular forms.

In this section we recall the formulas for the 
dimensions of some spaces of modular or cusp
forms associated to a representation $\rho$ of a discrete co-finite
subgroup of $\widetilde{SL}_2(\R)$. 
For weight at least 2 
the dimension is given by either the Riemann-Roch theorem or
the Selberg trace formula. More generally if $G$ is a group acting on
$A$ then it also acts on the spaces of cusp forms, and we calculate
the character of these representations.
These results are used in sections 9 and 12.

For weight $1/2$ forms Serre and Stark described an explicit basis 
as follows. 
\proclaim Theorem 7.1. 
Suppose that $\chi$ is an even Dirichlet character mod $N$. Then 
a basis for the space of modular forms of weight $1/2$ and character
$\chi_\theta\chi$ for $\Gamma_0(N)$ is given by the forms
$$\sum_{n\in \Z} \psi(n)q^{tn^2}$$ where $\psi$ is a primitive even
character of conductor $r(\psi)$, $t$ is a positive integer such that
$4r(\psi)^2t$ divides $N$, and $\chi(n)=\psi(n){D\choose n}$ for all
$n$ coprime to $N$, where $D$ is the discriminant of the quadratic
field $\Q[\sqrt{t}]$. (Note that $\psi$ is determined by $t$ and $\chi$.)

Proof. This is theorem A of [S-S, p.34]. 

The dimensions of spaces of holomorphic modular forms can all be
worked out as follows. For weight less than 0 there are no non-zero
forms, and weight 0 is trivial as these are just constants. For weight
$>2$ we can work out the dimension using the Selberg trace formula 
(see below) or
the Riemann-Roch theorem, and with a bit more care this also works for
weight 2 (there are extra correction terms coming from weight 0 forms
in this case). For weight $1/2$ the Serre-Stark theorem gives an
explicit basis, and this can be used to do the case of weight $3/2$
because the Selberg trace formula gives the difference of dimensions
for weights $k$ and $2-k$. This leaves the case of weight 1, which
seems to be the hardest case to do. In general weight 1 forms are
closely related to odd 2-dimensional complex representations of the
Galois group of $\Q$. Fortunately, for the low level cases we are
interested in, the weight 1 forms are usually easy to construct
explicitly using Eisenstein series and theta series of 2-dimensional
lattices (mainly because the exotic Galois representations only occur
for higher levels). 

Now we use the Selberg trace formula to find the dimensions of spaces
of forms of weight at least 2. 

If $X$ is a finite order automorphism of a finite dimensional complex
vector space $V$ with eigenvalues $\e(-\beta_j)$ for $1\le j\le \dim(V)$
and $0\le \beta_j<1$, then we define $\delta_\infty(X)$ to 
be $\sum(1/2-\beta_j)$, and we define $\delta_N(X)$ to be
$\delta_\infty(X)-\dim(V)/2N$. 
More generally, if $g$ is an endomorphism of $V$ commuting with 
the action of $G$ then we define $\delta_{\rho,\infty}(X,g)$
to be 
$$\sum (1/2-\beta_j)\Tr(g|V^{\e(\beta_j)X})$$
where the sum is over the distinct eigenvalues $\e(-\beta_j)$ of $X$, 
and we put $\delta_N(X,g)=
\delta_\infty(X,g)-\Tr(g)/2N$. 

\proclaim Lemma 7.2. If $\rho$ is a representation 
of a group containing $X$ on a finite dimensional 
complex vector space and $X^N=1$ then 
$$\eqalign{
\delta_N(X,g) 
&= {1\over N}\sum_{0<j<N}{\Tr_\rho(X^jg)\over 1-\e(j/N)}\cr
\delta_\infty(X,g) 
&= {\Tr(g)\over 2N} + {1\over N}\sum_{0<j<N}{\Tr_\rho(X^jg)\over 1-\e(j/N)}\cr
}$$

Proof. 
The trace of $g$ on 
the subspace of $\rho$ on which $X^{-1}$ has eigenvalue
$\e(k/N)$ is
$${1\over N}\sum_{j\bmod N}\Tr_\rho(gX^j)\e(jk/N).$$
Therefore 
$$\eqalign{
\delta_N(X)
&= \Tr_\rho(g)(1/2-1/2N)-
\sum_{0\le k<N}{k\over N}\times
 {1\over N}\sum_{j\bmod N}\Tr_\rho(gX^j)\e(jk/N)\cr
&= \Tr_\rho(g)(1/2-1/2N)-
{1\over N}\sum_{j\bmod N}\Tr_\rho(gX^j)\sum_{0\le k<N}{k\over N}\e(jk/N)\cr
&= 
{1\over N}\sum_{0<j<N}{\Tr_\rho(gX^j)\over 1-\e(j/N)}.\cr
}$$
This proves lemma 7.2. 

We will write 
$ModForm(\Gamma,k,\rho)$ for the space of modular forms of weight 
$k$ and representation $\rho$ for $\Gamma$ that are holomorphic
at cusps.

\proclaim Lemma 7.3. Suppose that $\Gamma$ is a discrete subgroup
of $\widetilde{SL}_2(\R)$ of co-finite volume and containing $Z$. 
Let $\rho$ be a complex representation of $\Gamma$ of finite dimension $d$ on
which $Z$ acts multiplication by some constant.
Choose $k\in \Q$ with $k>2$.
Then the dimension of the space of $ModForm(\Gamma, k, \rho)$
is equal to 0 unless $Z$ acts as $\e(- k/2)$, in which case it is 
$$(k-1)\dim(\rho)\omega(F)/4\pi
+\sum_{1\le j\le \rho}\delta_{\rho,\nu_j}(\e(k/2\nu_j)R_j)
+\sum_{1\le j\le \tau} \delta_{\rho,\infty}(T_j)
$$
where 
\item{$\omega(F)$} is the hyperbolic area of the fundamental domain 
$F$, and is equal to 
$$2\pi \Big( 2g-2 +\sum_{1\le j\le \rho} (1-1/\nu_j) +\tau \Big)$$
\item{$g$} is the genus of the compactification of $\Gamma\backslash H$
\item{$\rho$} is the number of elliptic fixed points in a fundamental 
domain.
\item{$\nu_j$} is the order of the $j$'th elliptic fixed point,
so the subgroup of $\Gamma$ fixing the $j$'th elliptic fixed point
is cyclic, generated by $R_j$ with $R_j^{\nu_j}=Z$. 
\item{$R_j$} is
the primitive elliptic element corresponding to $j$.
An elliptic element is called primitive
if it is conjugate to the ``clockwise'' element
$$
\left(
\pmatrix{ \cos(\pi/\nu)&-\sin(\pi/\nu)\cr
\sin(\pi/\nu)& \cos\pi/\nu}
, \sqrt{sin(\pi/\nu)\tau+\cos\pi/\nu}
\right)
\in \widetilde{SL}_2(\R)
$$ fixing $i$. 
\item{$\tau$} is the number of orbits of cusps of $\Gamma$ acting on $H$.
\item{$T_j$} is the unique element conjugate to $T^{-1}$ under 
$\widetilde{SL}_2(\R)$ such that the stabilizer of
the $j$'th cusp is generated by $T_j$ and $Z$. 

Proof. If $Z$ acts as multiplication by some constant not equal to 
$\e(-k/2)$ then the transformation of modular forms under
the element $Z$ immediately shows that any modular form
of weight $k$ is 0. 
If $Z$ acts as $\e(- k/2)$ then we get a ``multiplier system''
$\chi$ 
in the sense of [F, 1.3.4] from the representation $\rho$ by 
putting $\chi\left(\pmatrix{a&b\cr c&d\cr}\right)
 = \rho\left(\left(\pmatrix{a&b\cr c&d\cr}, \sqrt{c\tau+d}\right)\right)$
where we choose the value of the square root with $-\pi<\arg(\sqrt{c\tau+d})
\le \pi $. (Note that the variable $k$ used in [F] is half that used 
here.)

By [F, theorem 2.5.5] the dimension of the space of modular forms
is given by 
$$
(k-1)d\omega(F)/4\pi
+\sum_{1\le j\le \rho} 
\left( {d\over 2}-{d\over 2\nu_j} -{\alpha_j\over \nu_j} \right)
+{d\tau\over2} -\sum_{1\le j\le \tau} \beta_j
$$
where 
\item{$d$} $=\dim(\rho)$
\item{$\alpha_j$} is the sum of the numbers $\alpha_{jp}$ for
$1\le p\le d$, where the eigenvalues of $R_j$ are
$\e(-(k/2+\alpha_{jp})/\nu_j)$ and $\alpha_{jp}\in\{0,1,\ldots,\nu_j-1\}$. 
([F, p. 66--68])
\item{$\beta_j$} is the sum of the numbers $\beta_{jp}$ ($1\le p\le d$)
where the eigenvalues of $T_j^{-1}$ 
are $\e(\beta_{jp})$ and $0\le \beta_{jp}<1$. 

It is easy to check that 
$$
{d\over 2}-{d\over 2\nu_j} -{\alpha_j\over \nu_j}
=
\delta_{\rho,\nu_j}(\e(k/\nu_j)R_j)
$$
and
$$
{d\over 2} -\beta_j
=
\delta_{\rho,\infty}(T_j)
$$

Putting everything together proves lemma 7.3. 

Remark. If the weight is greater than 2 then
the dimension of the space of all cusp forms
is given by subtracting the dimension of the space of Eisenstein series, 
which is the sum over all cusps $j$ of the 
of the dimension of the subspace of $\rho$ fixed by $T_j$.

\proclaim Corollary 7.4. 
Suppose that $\rho$ is a finite dimensional representation of $\Gamma$.
Suppose that $k\in \Q$ and $k>2$. 
Then 
$$\dim(ModForm(\Gamma, k, \rho))=
{1\over 4}\sum_{0\le j<4} \e(k/2)\psi(Z^j)
$$
 where $\psi(g)$ is given by
$$
\psi(g)=
{(k-1)\omega(F)\over 4\pi} \Tr_\rho(g)
+\sum_{1\le j\le \rho} \delta_{\rho,\nu_j}(\e(k/2\nu_j)R_j,g)
+\sum_{1\le j\le \tau} \delta_{\rho,\infty}(T_j,g)
.
$$

Proof. Break up $\rho$ into the eigenspaces of $Z$ and apply 
lemma 7.3 to each eigenspace. This proves corollary 7.4. 

\proclaim Corollary 7.5. 
Suppose that $\rho$ is a finite dimensional representation of $\Gamma$
acted on by a finite group $G$.
Suppose that $k\in \Q$ and $k>2$. Then the character
of $ModForm(\Gamma, k, \rho)$, considered as a representation of $G$, 
is given by 
$$ \Tr(g|ModForm(\Gamma, k, \rho))
=
{1\over 4}\sum_{0\le j<4}\e(jk/2)\psi(gZ^j)
$$
where $\psi$ given by the formula of corollary 7.4 and $g\in G$. 

Proof. 
Note that the dimension of $ModForm(\Gamma, k, \sigma)$ is given 
by $\Tr(M|\sigma)$ for some element $M$ in the group ring of
$G$, whenever $\sigma$ is a representation of $\Gamma$ satisfying
the conditions of corollary 7.4. Therefore
$$\eqalign{
&\Tr(g|ModForm(\Gamma, k, \rho))\cr
=& \sum_{\chi\in Irred(G)}
\chi(g)\dim ((ModForm(\Gamma, k, \rho)\otimes \bar\chi)^G)\cr
=& \sum_{\chi\in Irred(G)}
\chi(g)\dim (ModForm(\Gamma, k, (\rho\otimes \bar\chi)^G))\cr
=& \sum_{\chi\in Irred(G)}
\chi(g)\Tr(M|(\rho\otimes \bar\chi)^G)\cr
=& \Tr(Mg| \rho).\cr
}$$
(The sums are over the sets $Irred$ of irreducible representations of $G$,
and $\bar\chi$ is the dual of the representation $\chi$.)
Therefore to find the trace of $g$, we just insert an extra factor of
$g$ whenever we have a trace in the formula for $\psi$.
This proves corollary 7.5. 

\proclaim 8.~The geometry of $\Gamma_0(N)$.

We summarize some standard results about the cusps and elliptic points
of the group $\Gamma_0(N)$. We need this information in order to use
the formulas of section 8.  For proofs see [Sh] or [Mi].

The group $\Gamma_0(N)$ has index $N\prod_{p|N}(1+1/p)$
(where the product is over all primes dividing $N$).

The equivalence class of the cusp $a/c$ of $\Gamma_0(N)$ is determined
by the invariants
$(c,N)$ (a divisor of $N$) and $(c/(c,N))^{-1}a$ (an element 
of $ (\Z/(c,N/c))^*$). A complete set of representatives for the cusps
is given by $a/c$ for $c|N$, $c>0$, $0<a\le (c,N/c)$, $(a,c)=1$. 
The cusp $a/c$ has width $N/(c^2,N)$.

In the rest of this section we work out the values of the characters
$\chi_\theta$ and $\chi_p$ on the elements $R_j$ and $T_{a/c}$ of section 7
associated to elliptic points or cusps of $\Gamma_0(N)$. 

\proclaim Lemma 8.1. If $(a,c)=1$ then the element $T_{a/c}$ of
$\widetilde{\Gamma}_0(N)$ is given by
$$T_{a/c} = 
\left(\pmatrix{1+act&-a^2t\cr c^2t&1-act\cr}, \sqrt{c^2t\tau+1-act}\right)
$$
where $t=N/(c^2,N)$. 

Proof. We can assume $c>0$ as the case $c<0$ follows from this, and the case
$c=0$ is trivial to check. 
The element $T_{a/c} $ is the conjugate of some element of the
form $\left(\pmatrix{1&-t\cr0&1\cr},1\right)$ for $t>0$, so is equal to 
$$
\left(\pmatrix{a&b\cr c&d\cr},\sqrt{c\tau+d}\right)
\left(\pmatrix{1&-t\cr 0&1\cr},1 \right)
\left(\pmatrix{d&-b\cr -c&a\cr},\sqrt{-c\tau+a}\right)
$$
where $t$ is the smallest positive integer such that this element is
in $\widetilde{\Gamma}_0(N)$, and $b$ and $d$ in $\Z$ are chosen so that
$ad-bc=1$. If we evaluate this element (keeping careful track of
the values of the square roots) we find it is equal to the expression in the 
lemma. This is in $\widetilde{\Gamma}_0(N)$ if and only if 
$N|c^2t$, so $t=N/(c^2,N)$. This proves lemma 8.1. 

\proclaim Lemma 8.2. If $p$ is an odd prime 
dividing $N$ then 
$\chi_p(T_{a/c})=1$. 

Proof. By lemma 8.1 we see that $\chi_p(T_{a/c})$ is 1 provided that
$1-act$ is a square mod $p$. But this is always true because $p$ must
divide either $c$ or $t=N/(c^2,N)$ if it divides $N$. This proves lemma 8.2. 

\proclaim Lemma 8.3. Suppose that $8|N$ and $(a,c)=1$. 
Then $\chi_2(T_{a/c})$ 
is $-1$ in the 3 cases 
$2||c$ and $8||N$, or $4||c$ and $8||N$, or $4||c$ and $16||N$, 
and is $1$ otherwise.

Proof. By lemma 8.1 we have 
$$\chi_2(T_{a/c})
=\chi_2
\left(\pmatrix{1+act&-a^2t\cr c^2t&1-act\cr}, \sqrt{c^2t\tau+1-act}\right)
={1-act\choose 2}
,$$
so
$\chi_2(T_{a/c})$ is 1 if 
$1-act\equiv 1\bmod 8$ and is $-1$ if $1-act\equiv 5\bmod 8$
(where $t=N/(c^2,N)$). 
Note that $act\equiv 1\bmod 4$ as $8|N$.
We do a case by case check to show that 
$act$ is $4\bmod 8$ in 
the 3 cases listed above, and is $0\bmod 8$ otherwise. 
\item{1.} If $c$ is odd then $8|t=N/(c^2,N)$, so we can assume that $c$ is even
(and hence $a$ is odd). 
\item{2.} If $8|c$ then $act\equiv 0\bmod 8$, so we can assume that
$2||c$ or $4||c$. 
\item{3.} If $2||c$ and $16|N$ then $4|t$ so $act\equiv 0\bmod 8$.
\item{4.} If $2||c$ and $8||N$ then $2||t$
so $act\equiv 4\bmod 8$.
\item{5.} If $4||c$ and $32|N$ then $2|t$ so $act\equiv 0\bmod 8$.
\item{6.} If $4||c$ and $32\not|N$ then $t$ is odd so $act\equiv 4\bmod 8$.

This proves lemma 8.3. 

\proclaim Lemma 8.4. Suppose that $4|N$. 
If $2||c$ and $4||N$ then $\chi_\theta(T_{a/c})=-i^{t}$, with 
$t=N/(c^2,N)$. 
If $2||c$ and $8||N$ then $\chi_\theta(T_{a/c})=-1$.
Otherwise $\chi_\theta(T_{a/c})=1$ 

Proof. By lemmas 8.1 and 5.2, we see that if $4|act$ then
$$\chi_\theta(T_{a/c})
=\chi_\theta
\left(\pmatrix{1+act&-a^2t\cr c^2t&1-act\cr}, \sqrt{c^2t\tau+1-act}\right)
={c^2t\choose 1-act}
={1-act\choose t}
.$$
This is equal to 1 unless $2||t$ and $2||c$, in which case it is $-1$
and $8||N$. So lemma 8.4 is true whenever $4|act$. 

If $c$ is odd then $4|t$ so $4|act$, and if $4|c$ then $4|act$, so we
can assume that $2||c$. If $8|N$ and $2||c$ then $2|t$ so
$4|act$. So we can also assume that $4||N$, which implies that $a$ and
$t$ are odd and $act\equiv 2\bmod 4$.
Then by lemmas 8.1 and 5.2, 
$$
\chi_\theta(T_{a/c})
=-i{t\choose 1-act}
=-i(-1)^{(t-1)(1-act-1)/4}{1-act\choose t}
=-i(-1)^{(t-1)/2}
=-i^t
.
$$

This proves lemma 8.4. 

Finally we work out the values of characters on elliptic elements. 
Note that if the characters $\chi_\theta$ or $\chi_2$ are non-trivial,
then $4|N$ so there are no elliptic elements. So we only have
to do the characters $\chi_p$ on elliptic elements for $p$ an odd prime. 

\proclaim Lemma 8.5. Suppose that $\Gamma_0(N)$ has an elliptic fixed point of
order 2 fixed by the primitive elliptic element 
$R\in \Gamma_0(N)$ as in lemma 7.3, 
 and let $p$ be an odd prime dividing $N$
(so that $p\equiv 1\bmod 4$). 
Then $\chi_p(R) =(-1)^{(p-1)/4}$. 

Proof. We know that $R^2=\pm \pmatrix{1&0\cr 0&1\cr}$ 
so $R=\pmatrix{a&b\cr c&-a}$
for some $a,b,c$, and $\chi_p(R)={a\choose p}$. As $-a^2-bc=1$ and $N|c$
and $p|N$
we see that $a^2\equiv -1\bmod p$, so that $a$ is an element of order 
exactly 4
in $(\Z/p\Z)^*$, and hence is in the (unique) index 2 subgroup 
of $(\Z/p\Z)^*$ if and only if $p\equiv 1\bmod 8$. But 
$\chi_p(a)=1$ if and only if $a$ is in the index 2 subgroup of $(\Z/p\Z)^*$. 
This proves lemma 8.5. 

\proclaim Lemma 8.6. Suppose that $\Gamma_0(N)$ has an elliptic fixed point of
order 3 fixed by the primitive elliptic element 
$R\in \Gamma_0(N)$ as in lemma 7.3, 
and let $p$ be an odd prime dividing $N$. 
Then $\chi_p(R) =(-1)^{(p-1)/2}$. 

Proof. The character $\chi_p$ has order 2, so 
$$\chi_p(R)=\chi_p(R^3)=\chi_p(Z) = {-1\choose p}= (-1)^{(p-1)/2}.$$
This proves lemma 8.6. 

\proclaim 9.~An application of Serre duality.

In this section we summarize some results from [B99]. 
We will later often need to find modular forms with given singularities
at cusps. In this section we show how it is sometimes possible to
prove the existence of modular forms with given singularities without
writing them down explicitly. The main idea is that Serre duality 
shows that the space of obstructions to finding a modular form with given 
singularities at cusps is a space of cusp forms whose dimension can usually
be worked out explicitly. If this space of obstructions has smaller dimension
than the space of potential solutions, then at least one solution must exist.

If $\kappa$ is a cusp of $\Gamma$, let $q_\kappa$ be a uniformizing parameter
at $\kappa$ on $\Gamma\backslash H$. For a representation $\rho$ on
$V_\rho$,
let $V_{\rho}^*$ denote the dual.
Let
$$PowSer_\kappa(\Gamma,\rho)=\C[[q_\kappa]]\otimes V_\rho$$
be the space of formal power series in $q_\kappa$
with coefficients in $V_\rho$, 
let
$$Laur_\kappa(\Gamma,\rho)=\C[[q_\kappa]][q_\kappa^{-1}]\otimes V_\rho$$
be the space of formal Laurent series in $q_\kappa$
with coefficients in $V_\rho$, 
and let
$$
Sing_\kappa(\Gamma,\rho) 
= {Laur_\kappa(\Gamma,\rho)\over PowSer_\kappa(\Gamma,\rho)}
$$
be the space of possible 
singularities 
of $V_\rho$ valued Laurent series at $\kappa$. 
The two
spaces $PowSer_\kappa(\Gamma,\rho^*)$ 
and $Sing_\kappa(\Gamma,\rho)$ are paired into $\C$ by taking the residue
$$\langle f,\phi\rangle = 
\hbox{Res}\big( f \phi\, q_\kappa^{-1} \,dq_\kappa \big),$$
for $f\in PowSer_\kappa(\Gamma,\rho^*)$ and $\phi\in Sing_\kappa(\Gamma,\rho)$.
Here the product of $f$ and $\phi$ 
is defined using the pairing of $V_\rho$ and $V_\rho^*$.

Then the spaces
$$Sing(\Gamma,\rho) = \oplus_{\kappa} Sing_\kappa(\Gamma,\rho)$$
and
$$PowSer(\Gamma,\rho^*) = \oplus_{\kappa} PowSer_\kappa(\Gamma,\rho^*),$$
where $\kappa$ runs over the $\Gamma$-inequivalent cusps, are paired by the
sum of the local pairings at the cusps.

There are maps
$$\lambda:CuspForm(\Gamma,k,\rho^*) \longrightarrow PowSer(\Gamma,\rho^*)$$
and
$$\lambda:SingModForm(\Gamma,2-k,\rho) \longrightarrow Sing(\Gamma,\rho),$$
defined in the obvious way
by taking the Fourier expansions of their nonpositive part at the
various cusps.

We define the space $Obstruct(\Gamma,k,\rho)$
of obstructions to finding a modular form of type
$\rho$ and weight $k$ which is holomorphic on $H$ and has given
meromorphic singularities at the cusps to be the
space
$$Obstruct(\Gamma,k,\rho) =
{Sing(\Gamma,\rho)
\over \lambda(SingModForm(\Gamma,k,\rho))}.
$$

\proclaim Lemma 9.1. Suppose that $\rho$ is a representation 
of $\Gamma$ factoring through some finite
quotient of $\Gamma$ and $k\in \Q$ with $Z=\e(-k/2)$ on $\rho$. Then
$$Obstruct(\Gamma, 2-k, \rho)$$ 
is dual to 
$$CuspForm(\Gamma, k, \rho^*)
$$
(and both spaces are finite dimensional). 

Proof. This can be proved in the same way as theorem 3.1 of [B99]
(which is really a special case of Serre duality). The only real
difference is that in [B99] the space $Sing$ is different because it
also includes information about the constant terms of functions. This
has the effect of replacing the space of holomorphic modular forms in
[B99] by a space $CuspForm(\Gamma, k, \rho^*)$ of cusp forms. This
proves lemma 9.1.

{\bf Example 9.2.} In one of the examples later on, we need to know
that there is a non-zero weight $-7$ form of character $\chi_3$ for
$\Gamma_0(3)$ whose singularities are a pole of order at most 1 at
$i\infty$ and a pole of order at most 3 at 0, and such that the
coefficient of $q_3^2$ at the cusp 0 vanishes. The space of possible
singularities is spanned by $q^{-1}$ at $i\infty$ and by $q_3^{-1}$
and $q_3^{-3}$ at 0, so it is 3 dimensional. The space of obstructions
is the space of weight $2--7=9$ cusp forms of character $\chi_3$ for
$\Gamma_0(3)$, which has dimension 2. This is less than 3, so the
space of forms with the singularities above is at least one
dimensional, so a nonzero form exists. (For an explicit formula for
it see section 12.)

Warning. This method only gives a lower bound for dimensions of spaces
of forms with given singularities. In practice this lower bound is 
often the exact dimension, but there are occasional examples 
where the lower bound is 0 but nevertheless there is a nonzero form. 

\proclaim 10.~Eisenstein series.

We summarize some standard results about Eisenstein series that we will use
in section 12. 

\proclaim Lemma 10.1. Assume that $k\ge2$ is even.
Put
$$E_k(\tau) = 1-{2k\over B_{k}}\sum_{n\ge 1} q^n\sum_{d|n} d^{k-1}.$$
If $k>2$ then $E_k$ is a modular form of weight $k$ for $\Gamma_0(1)$. 
If $\sum_{d|N}a_d/d=0$ then $\sum_{d|N} a_dE_2(\tau)$ is a modular form
of weight 2 for $\Gamma_0(N)$. 

Proof. These are just the usual Eisenstein series for $SL_2(\Z)$. 

\proclaim Lemma 10.2. Assume that $k\ge2$ is integral and let $\chi$
be a non-principal Dirichlet character mod $N$ with $\chi(-1)=(-1)^k$. Then 
$$E_k(\tau,\chi) = \sum_{n\ge 1} q^n\sum_{d|n} d^{k-1}\chi(n/d)$$
is a modular form of weight $k$ and character $\chi$ for $\Gamma_0(N)$. 

Proof. See [Mi, theorems 7.1.3 and 7.2.12 and lemma 7.1.1]. 

\proclaim Lemma 10.3. Let $\chi$
be a non-principal Dirichlet character mod $N$ with $\chi(-1)=-1$. Then 
$$E_1(\tau,\chi) =1+{2\over L(0,\chi)} \sum_{n\ge 1} q^n\sum_{d|n}\chi(n/d),$$
where
$$L(0,\chi)= -B_{1,\chi}=-\sum_{0<n<N} \chi(n)n/ N,$$
is a modular form of weight $k$ and character $\chi$ for $\Gamma_0(N)$. 

Proof. See [Mi, theorem 7.2.13 and lemma 7.1.1]. 

{\bf Example 10.4.} 
The Eisenstein series $E_1(\tau,\chi_3)$ is a weight 1 modular form 
for $\Gamma_0(3)$ of character $\chi_3$ whose power series is
$$
E_1(\tau,\chi_3)
=1+6\sum_{n>0}q^n\sum_{d|n}\chi_3(n/d)
=1+6q+6q^3+6q^4+O(q^7).
$$
This form is also the theta function of the $A_2$ lattice. (It is  used
in Wiles' proof of Fermat's last theorem to show that
any weight 1 form is congruent mod 3 to a weight 2 form: just multiply
by $E_1(\tau,\chi_3)$!)

\proclaim 11.~Reflective forms.

Suppose that $L$ is a lattice of level dividing $N$. 
We define a singularity at a cusp $a/c$ of $\Gamma_0(N)$ to be {\bf reflective}
if it is a linear combination of terms $q^{-1/n}$ for
$n$ a positive integer, such that
every norm $-2/n$ element of $(L'/L)^{c*}$ has order dividing $n$. 
Note that if a norm $-2/n$ element of $L'/L$ has order dividing
$n$ then every lift of it to a norm $-2/n$ element of $L'$ is a root
of $L'$. 
We say that a modular form is {\bf reflective} for $L$
if it is a meromorphic modular form of weight $\sign(L)/2$, level
$N$, character $\chi_L$, and its only singularities
are reflective singularities at cusps. 

This definition is the result of a lot of experimentation to find a
definition that is easy to use and that also seems to give most of the
``interesting'' lattices. There are many possible variations of it,
some of which we now briefly describe. First of all, we could use
vector valued modular forms rather than scalar valued modular forms.
The main problem with this is practical inconvenience: it is not that
easy (though possible) to work with modular forms taking values in a
space of dimension (say) $2^9$. Moreover examples suggest that the
most interesting vector valued modular forms are often invariant under
$\Aut(A)$. This suggests using invariant vector valued modular forms
instead, but examples suggest that these are very closely related to
scalar valued modular forms of level $N$, so we may as well use these.
The ``allowable'' singularities can also be varied. For example, we
could also allow singularities of the form $q^{-n}$ such that $A$ has
no vectors of norm $-2n$. One problem with this is that it turns out
there are then too many modular forms with these singularities if the
$p$ rank of $A$ is 1 for some prime $p$. Another possibility is to
allow singularities of the automorphic form that correspond to
characteristic vectors as well as roots. The reason for this is that
this would include the reflection groups of the lattices $I_{1,20}$,
$I_{1,21}$, $I_{1,22}$, $I_{1,23}$ that were described in [B87].

The main theme of this paper is that
lattices of negative signature that have a non-zero reflective 
modular form tend to be ``interesting''. The meaning of ``interesting''
depends on the dimension  of the maximal positive definite sublattice. 
For example, negative definite lattices might be similar to the Leech lattice, 
Lorentzian lattices should have interesting reflection groups, and lattices 
in $\R^{2,n}$ should have interesting automorphic forms associated with them.
The main reason for 
the definition of reflective forms is the following lemma. 

\proclaim Lemma 11.1. Suppose that $f$ is a reflective modular form for
the lattice $L$. Then all singularities of the automorphic form 
of $f$ on $G(L)$ are orthogonal to negative norm roots of $L$. 

Proof. By theorem 6.2 of [B98a], the singularities of $\Phi$
are orthogonal to vectors $\gamma\in L'$ such that 
$\gamma^2<0$ and the vector valued form 
$$\sum_{g\in \widetilde{SL}_2(\Z)/\Gamma_0(N)} g(f)\rho_A(g)(1)$$
of $f$ has a singularity
of type $q^{\gamma^2/2}e_\gamma$. 
By lemma 3.2, if the vector valued form has such a singularity, 
then $f$ has a singularity $q^{\gamma^2/2}$ at some cusp $a/c$ with 
the image of $\gamma\in A^{c*}$. 
By the definition of a reflective form this implies that $\gamma$ is a
root of $L'$. 
This proves lemma 11.1.

Warning. If a level $N$ modular form $f$ is nonzero it is possible
that the corresponding vector valued modular form is zero.  And even
if the vector valued modular form has nontrivial singularities at
cusps, it is possible that the corresponding automorphic form is zero.
Although we normally expect a sing of a vector valued form to give a
sing of the corresponding automorphic form, it is possible that all
such singularities correspond to vectors of an empty set, or it is
possible that all singularities happen to cancel each other out.
(This is related to the well known problem that the image of a
non-zero form under the theta correspondence can be zero.)

For lattices of some given level $N$ it is usually easy to determine
the reflective singularities, though the answer sometimes involves
many different cases. The following lemma gives a simple condition for
a singularity to be reflective, which in practice covers many cases.

\proclaim Lemma 11.2. Suppose $a/c$ is a cusp of $\Gamma_0(N)$ and $L$
is an even lattice of level equal to $N$. 
Then a pole of order 1 at $a/c$ is reflective if $(c,N)$ is a Hall divisor
of $N$. 

Proof. The cusp $a/c$ has width $h=N/(c^2,N)= N/(c,N)$. 
To show that $q_h^{-1}= q^{-1/h}$ is a reflective singularity at $a/c$
it is enough to show that all elements $\alpha\in A^{c*}$ with
$(\alpha,\alpha)\equiv -2/h\bmod 2$ satisfy $h\alpha=0$. But this
is obvious because $A^{c*}=A^c$ is the set of elements of
$A$ of order dividing $h$ because $(N,c)$ is a Hall divisor of $N$. 
This proves lemma 11.2. 

\proclaim 12.~Examples.

In this section we give some examples of lattices with reflective modular 
forms. For a given level the method for finding such forms is as follows:
\item{1.} 
Work out the ring of modular forms for $\Gamma_0^2(N)$, paying
particular attention to the forms with zeros only at cusps.
\item{2.} 
Work out the possible reflective singularities for each possible discriminant 
form of level $N$. 
\item{3.} 
Try to find modular forms with reflective singularities for each
discriminant form.
\item{4.} 
Try to find lattices corresponding to these modular forms.
\item{5.} 
For each lattice with a non-zero reflective modular form, 
see if it is connected with any interesting reflection groups, 
automorphic forms, moduli spaces, or Lie algebras. 

The examples probably include most interesting cases for small prime
level, but become less and less complete as the level gets larger
because the number of cases to consider becomes rather large. What
usually seems to happen is that for each level there is some
``critical'' signature, with the property that almost all lattices up
to that signature have non-zero reflective modular forms, but beyond
that signature there are only a few isolated examples, usually with
$p$-rank at most 2 for some prime $p$. For example, for level $N=1$
the critical signature is $-12$, corresponding to the Leech lattice,
the lattice $II_{1,25}$ whose reflection group was described by
Conway, and so on, while for level $N=2$ the critical signature is
$-16$ corresponding to the Barnes-Wall lattice and so on. 

There are several other methods for finding examples of lattices with
reflective modular forms. First, the inverses of eta quotients of elements of
Conway's group of automorphisms of the Leech lattice are often
reflective modular forms for various lattices. Second, Hauptmoduls
of genus zero subgroups of $SL_2(\R)$ are often reflective modular
forms for lattices of signature 0; see the case $N=17$ below for some
examples of this. More generally, several other modular functions of
genus 0 subgroups with poles of order 1 at some cusps are reflective
modular forms for some lattices. (The restriction to genus 0 subgroups
is just an observation: most cases seem to be related to genus 0
subgroups. I do not know of a good theoretical reason for this.)
Third, Y. Martin [M] gave a list of many eta quotients with
multiplicative coefficients, and again many of these seem to be the
inverses of reflective modular forms. (Note that there are many
multiplicative eta quotients not on Y. Martin's list that also appear as
reflective modular forms, because Y. Martin restricted himself to forms
of integral weight whose conjugate under the Fricke involution was
also multiplicative.) Y. Martin's list contains many examples high level,
up to level 576, which suggests that there should be many reflective
modular forms beyond the examples in this section.

Many of the calculations with modular forms in this section were done
using the PARI calculator [BBCO].

In most of the examples below, the tables have the following meaning. 
``Group'' describes the relevant subgroup of $SL_2(\Z)$, ``index''
gives its index in $SL_2{\Z}$, $\eta_2$, $\eta_3$, and $\eta_\infty$
are the numbers of elliptic points of orders 2 and 3 and the numbers of cusps
of a fundamental domain, and ``genus'' gives the genus of the corresponding 
compact Riemann surface. The second table in each section is a table of
the cusps, with one line per cusp. The column ``cusps'' gives a 
representative cusp, ``width'' is the width of the cusp, ``characters''
lists the non-trivial values of characters on the normalized generator of
the subgroup fixing a cusp, ``$\eta$'' lists an eta function 
with a zero of order ``zero'' at this cusp and no other zero and with
weight ``weight'' and character ``character''. 

The Hilbert function is the rational function of $x$ whose coefficients
give the dimensions of the spaces of modular forms of various weights. 
Sometimes we put in extra variables $u_{p}$ and $u_{\theta}$, which 
describe the dimensions of spaces with various characters.

We sometimes write $q_n$ for $\e(\tau/n)$.

\proclaim $N=1$.

All the results we get for this case are well known, but this case is
still useful as a warming up exercise.

The discriminant form $A$ has order 1, and the group $\tilde \Gamma_0(1)$
is just $\widetilde{SL}_2(\Z)$. The character $\chi$ is always 1. 

\vbox{
\halign{
#~&~#~&~#~&~#~&~#~&~#~\cr
Group         &index&$\nu_2$&$\nu_3$&$\nu_\infty$&genus\cr
$\Gamma_0( 1)$&    1&      1&      1&           1&    0\cr
}}
\halign{
~#~&~#~&~$#$~&~#~&~#~&~$#$~\cr
cusps     & width & \eta            & zero & weight  \cr
$i\infty$ &    1  &  1^{24}         &  1   & 12      \cr
}

The ring of modular forms is a polynomial ring with generators
given by the Eisenstein series
$E_4$ and $E_6$ of weights 4 and 6, and the Hilbert function is
$1/(1-x^4)(1-x^6)$.
The critical weight is $k=12$, and the critical form is
$\Delta(\tau)=\eta(\tau)^{24}$.
By looking at the forms $E_4(\tau)^n/\Delta(\tau)$ we see that
every even lattice $L$ with $N=1$ and $\sign(L)\ge -24$ has
a nonzero reflective modular form.
The signature must be $0\bmod 8$. 

Next we find some possible reflective forms.
By lemma 11.2, poles of order at most 1 at the cusp are reflective
singularities. The modular forms with poles of order at most
1 at the cusp are exactly those of the form
$(\hbox{holomorphic modular form})/\Delta$. So there are no examples
of weight less than $-12$, and the ones of weights $-12$, $-8$, and $-4$
are multiples of $1/\Delta$, $E_4/\Delta$, and $E_4^2/\Delta$.

We now look at some of these
cases in more detail.

For signature $-24$ we take $f$ to be
$1/\Delta(\tau)=q^{-1}+24+\cdots$.  We get an automorphic form for the
lattice $II_{2,26}$ of singular weight $24/2=12$.  The corresponding
Lie algebra is the fake monster Lie algebra. Its Weyl group is
Conway's reflection group of the lattice $II_{1,25}$, which has a norm
0 Weyl vector. The critical lattice of this Weyl vector is of course
the Leech lattice.

The lattice $II_{4,28}$ has a reflective form, and is the underlying
integral lattice of Allcock's largest quaternionic reflection group
[A].  The lattices $II_{n,16+n}$ for various values of $n>2$ appear in
the moduli space of K3 surfaces, possibly with some extra structure
such as a B-field. The existence of a reflective form for these
lattices appears to be significant in the corresponding moduli
spaces as it is usually necessary to discard points of the Grassmannian
that are orthogonal to norm $-2$ vectors.

For signature $-8$ and $-16$ we find the lattices
$II_{1,9}$ and $II_{1,17}$ whose (arithmetic) reflection groups
were first described by Vinberg [V75].

For signature $0$ we take $f$ to be $j(\tau)-744 = q^{-1}+196884q+\cdots$.
We get an automorphic function for the lattice $II_{2,2}$,
which is more or less $j(\sigma)-j(\tau)$ in suitable coordinates.
The corresponding Lie algebra
is the  monster Lie algebra. The reflection group is not very interesting
as it is of order 2 (and is the Weyl group of the monster Lie algebra).

\proclaim  $N=2$.

The discriminant form $A$ has order $2^{2n}$ for some
non-negative integer $n$. The character $\chi_A$ is always
trivial as $A$ always has square order and signature  divisible by 4.

The group $\Gamma_0(2)=\Gamma^2_0(2)$
has 2 cusps which can be taken as $i\infty$ (of width 1) and 0 (of width 2).
It has one
elliptic point of order 2, which  can be take as the point $(1+i)/2$,
fixed by $\pmatrix{1&-1\cr2&-1\cr}$.

\vbox{\halign{
#~&~#~&~#~&~#~&~#~&~#~\cr
Group         &index&$\nu_2$&$\nu_3$&$\nu_\infty$&genus\cr
$\Gamma_0( 2)$&    3&      1&      0&           2&    0\cr
}}
\halign{
~#~&~#~&~$#$~&~#~&~#~&~$#$~\cr
cusps     & width & \eta            & zero &  weight \cr
        0 &    2  &  1^{16}2^{-8}   &  1   &  4      \cr
$i\infty$ &    1  &  1^{-8}2^{16}   &  1   &  4      \cr
}

The ring of modular forms for $\Gamma_0^2(2)=\Gamma_0(2)$ is a polynomial ring
on the generators
$-E_2(\tau)+2E_2(2\tau)= \theta_{D_4}(\tau)= 1+24q+24q^2+O(q^3)$ of weight 2
with a zero at the elliptic point,
and $E_4(\tau)$ of weight 4.
The Hilbert function is $1/(1-x^2)(1-x^4)$. 

All poles of order at most 1 at cusps are reflective by lemma 11.2 as
$N=2$ is square-free. If $A=II(2^{-2})$ then $A$ has no non-zero elements
of norm $1\bmod 2$, so a pole of order 2 at the cusp $0$ is also reflective.
(There are also other possible reflective singularities for
lattices of high 2-rank.)

By looking at the form $\Delta_{2+}(\tau)^{-1}$, with order 1 poles at
all cusps, we see that all level 2 even lattices of signature at least
$-16$ have reflective modular forms.  The Lorentzian lattices
$II_{1,17}(2^{+8})$ and $II_{1,17}(2^{+10})$ have norm 0 Weyl vectors;
their reflection groups are not arithmetic, but are similar to the case of
$II_{1,25}$.
The remaining Lorentzian lattices of dimension at most 18 have
positive norm Weyl vectors, so their reflection groups are arithmetic.
The lattice $II_{1,17}(2^{+2})$ is the even sublattice of
an odd unimodular lattice; see [V75]. The
arithmetic reflection group of the lattice  $II_{1,17}(2^{+4})$
appeared recently in Kondo and Keum's work [K-K] as the Picard lattice
of the Kummer surface of a generic product of elliptic curves,
and can be obtained as the orthogonal complement of a $D_4^2$ in
$II_{1,25}$.

The reflection group of $II_{1,17}(2^{+6})$ seems to be the highest
dimension of a ``new'' example of an arithmetic reflection group in
this paper. This lattice has an usually complicated fundamental
domain, with 896+64 sides. It can be described as follows. Let $K$ be the
16 dimensional even lattice in the genus $II_{0,16}(2^{+6})$ that has
root system $A_1^{16}$. It can be constructed by applying construction
A of [C-S, chapter 5] to the first order Reed-Muller code of length 16
with 32 elements (see [C-S, p. 129], which uses construction B applied
to this code to construct the Barnes-Wall lattice). Then
$L=II_{1,17}(2^{+6})$ can be constructed as $K\oplus II_{1,1}$. The
fundamental domain of the reflection group $R$ of $L$ has 2 norm 0
vectors $z$, $z'$ of type $K$ (together with many other norm 0 vectors
of other lattices).  The generalized Weyl vector is given by
$\rho=(z+z')/2$. The fundamental domain domain has 64 walls
corresponding to norm $-2$ roots of $L$ and 896 walls corresponding
to norm $-4$ roots of $L$ (or equivalently to norm $-1$ roots of
$L'$). The norm $-2$ simple roots split into 2 groups of 32. The first
group of 32 have inner product $0$ with $z$ and $-1$ with $z'$ and
correspond to the 16 coordinate vectors of $K$ and their
negatives. The other 32 have inner product $-1$ with $z$ and $0$ with
$z'$ and correspond to the 32 elements of the Reed-Muller code. 
(There is of course an automorphism  of the fundamental domain exchanging
$z$ and $z'$ and the two groups of 32 norm $-2$ simple roots.) The
896 norm $-1$ simple roots of $L'$ all have inner product 1 with both
$z$ and $z'$ and correspond to the  elements of the
dual of the Reed-Muller code whose weight is $2\bmod 4$. The Reed-Muller code
has Hamming weight enumerator $x^{16}+30x^8y^8+y^{16}$, so by
the MacWilliams identity ([C-S, p.78]) the dual code has weight enumerator
$$
\eqalign{
&{(x+y)^{16}+30(x+y)^8(x-y)^8+(x-y)^{16}\over32} \cr
=
& x^{16} +
140y^4x^{12} + 448y^6x^{10} + 870y^8x^8 + 448y^{10}x^6 + 140y^{12}x^4
+ y^{16}\cr
}$$ and therefore there are 448+448=896 elements of length
$2\bmod 4$.  Alternatively, $L$ can be constructed as the orthogonal
complement of a certain $A_1^8$ in the Dynkin diagram of $II_{1,25}$.
(Note that the Dynkin diagram of $II_{1,25}$ contains more than
one orbit of subsets isomorphic to $A_1^8$. The orbit 
we use has the special property that it is not contained in 
a $A_1^7A_2$ sub-diagram.)  The 64 norm $-2$ roots
correspond to the 64 ways to extend the $A_1^8$ to an $A_1^9$, and the
896 norm $-1$ simple roots correspond to the 896 ways of extending it
to an $A_3A_1^6$ diagram.  These $A_3A_1^6$ diagrams are the same as
those used by Kondo in [K]  to describe the automorphism group of a
generic Jacobian Kummer surface.  The automorphism group of the Dynkin
diagram of $L$ has a normal subgroup of order $2^{10}$ and the
quotient is the alternating group $A_8$.
Unfortunately 
$II_{1,17}(2^{+6})$ cannot be the Picard
lattice of a K3 surface:  S. Kondo pointed out to me that
its 2-rank (6) is larger than
its codimension (4) in $II_{3,19}$.

The lattices $II_{n,n+20}(2^{-2})$ also have a reflective modular form
$\theta_{D_4}(\tau)/\Delta(\tau)$.  In particular the Lorentzian
lattice $II_{1,21}(2^{-2})$ has an arithmetic reflection group; see
[B87, p.149] for a description of its fundamental domain. Esselmann [E]
showed that this is essentially the only example of a cofinite
reflection group of a Lorentzian lattice of dimension at least 21.
(Of course we can find trivial variations of this example by
taking the Atkin-Lehner conjugate $II_{1,21}(2^{-20})$, or by multiplying
all norms by a constant.)

We can find some automorphic forms of singular weight, corresponding
to lattices $L$ and reflective modular forms $f$ as follows.
\item{1} $L=II_{2,18}(2^{+10})$, $f=\eta_{1^{-8}2^{-8}}$.
This is the denominator function of a generalized Kac-Moody algebra of
rank 18. This example is related to the element $2A$ of
$\Aut(\Lambda)$, of cycle shape $1^82^{8}$.  There are 24 lattices in
the genus $II_{0,16}(2^{+8})$ by [S-V], one of which is the Barnes-Wall
lattice, and the others all have root systems of rank 16.
\item{2} $L=II_{2,10}(2^{+2})$, $f=2^4\eta_{1^{-16}2^8}$.
This example is related to the element $-2A$ of $\Aut(\Lambda)$.
\item{3} $L=II_{2,10}(2^{+10})$, $f=\eta_{1^{8}2^{-16}}$

The first case is closely related to the reflection group of the
lattice $II_{1,17}(2^{+8})$, whose fundamental domain has a nonzero
norm 0 vector fixed by its automorphism group, as in the lattice
$II_{1,25}$. For more about this lattice and its reflection group see
[B90].

The last two cases are really the same, since the lattices are
Atkin-Lehner conjugates of each other, and the automorphic forms we
get are more or less the same. This automorphic form is the
denominator function of two generalized Kac-Moody superalgebras,
and is also closely related to the moduli space of Enriques surfaces.
See [B98a, example 13.7] and [B96] for more details.

If $R$ is the reflection group of the lattice $II_{1,17}(2^{+2})$
generated by the reflections of norm $-2$ vectors and $D$ is its
fundamental domain, then $\Aut(D)$ has a finite index subgroup
isomorphic to $\Z$ and fixes a nonzero norm 0 vector $z$. However
there seems to be no reflective modular form for $\Gamma_0(2)$
corresponding to this reflection group.

The remaining level 2 cases of signatures $-4$, $-8$ and $-12$
are left to the reader; they  all give known arithmetic reflection groups,
often associated to unimodular lattices as in [V75]. 

\proclaim $N=3$.

\vbox{\halign{
#~&~#~&~#~&~#~&~#~&~#~\cr
Group         &index&$\nu_2$&$\nu_3$&$\nu_\infty$&genus\cr
$\Gamma_0( 3)$&    4&      0&      1&           2&    0\cr
}}
\halign{
~#~&~#~&~$#$~&~#~&~#~&~$#$~\cr
cusps     & width & \eta            & zero & weight \cr
        0 &    3  &  1^{ 9}3^{-3}   &  1   & 3      \cr
$i\infty$ &    1  &  1^{-3}3^{ 9}   &  1   & 3      \cr
}
The character $\chi_3$ is trivial for forms of even weight
and nontrivial for forms of odd weight. The forms of integer
weights and arbitrary character are the same as the forms for
$\Gamma_1(3)$ and trivial character. Note that
$\Gamma_0(3)$ is the product of $\Gamma_1(3)$ and its center of
order 2 generated by $Z$.

The ring of modular forms for $\Gamma_0^2(3)= \Gamma_1(3)$ is
a polynomial ring generated by
$\theta_{A_2}(\tau)= E_1(\tau, \chi_3) = 1+6q+6q^3+6q^4+O(q^7)$ of weight 1 and
$E_3(\tau,\chi_3)=\eta(\tau)^{-3}\eta(3\tau)^9= 
q + 3q^2 + 9q^3 + 13q^4 + 24q^5 + O(q^6)$
of weight 3.
The Hilbert function is
$1/(1-u_3x)(1-u_3x^3)$.

Next we find some reflective singularities. At the cusp $i\infty$
the singularity $q^{-1}$ is reflective. At the cusp $0$ the singularity
$q_3^{-1}$ is reflective.
If the discriminant form $A$ is $II(3^{\pm 1})$ or   $II(3^{-2})$
then $A$ has  no nonzero elements
of norm $0\bmod 2$,  so the singularities $q_3^{-1}$, $q_3^{-3}$
are reflective at the cusp 0.

The forms $\Theta_{A_2}(\tau)^n/\Delta_{3+}$ show  that all level 3 even
lattices of signature at least $-12$ have non-zero reflective modular forms.

There are also some other examples of reflective modular forms
for lattices of small 3-rank.
If we take the signature to be $-18$ and take $A$ to be $II(3^{+1})$
then the form $\Theta_{E_6}(\tau)/\Delta(\tau)$ is reflective.
This can be used to show
that the reflection group of the lattice $II_{1,19}(3^{+1})$
is arithmetic. This reflection group was first found by Vinberg [V85]
in his investigations of the ``most algebraic'' K3 surfaces.
The lattice can also be constructed as the orthogonal complement of an $E_6$
in $II_{1,25}$, and this gives another proof that the
reflection group is arithmetic [B87]. 

Next take the signature to be $-14$ and take $A$ to be $II(3^{-1})$.
We let $f$ be the form
$$
{
E_1(\tau,\chi_3)^5-270E_1(\tau, \chi_3)^2\eta(\tau)^{-3}\eta(3\tau)^9
\over
\Delta(\tau)
}
=q^{-1} -216 -9126q+O(q^2).
$$
The constant 270 is chosen so that the coefficient of $q^{-2/3}=q_3^{-2}$ of
$f(-1/\tau)=-i3^{-5/2}\tau^3 (-9q^{-1}+810q^{-1/3}+ 1944 + 53136q^{2/3}+O(q))$
vanishes. So  the automorphic forms with singularities
constructed from $f$ have all their singularities orthogonal to roots.
However something unexpected now happens: the $\C[A]$
valued modular form induced from $f$ is identically zero! So the
piecewise linear automorphic forms  constructed from $f$
as in theorem 4.2 have no singularities and are also
zero.

In spite of this the reflection group of $L=II_{1,15}(3^{-1})$ still
has a nonzero vector (of norm 0) in the fundamental domain fixed by the
automorphism group of the fundamental domain. To see this, we
represent $L$ as the orthogonal complement of an $A_2$ in
$II_{1,17}$. Then the quotient of $\Aut^+(L)$ by the reflection group
can be worked out using theorem 2.7 of [B98b] and turns out to be an
infinite dihedral group, which has an index 2 subgroup isomorphic to
$\Z$.  Next we can classify the primitive norm 0 vectors $z$ of $L$,
and we find that there are just two orbits, with the lattice
$z^\perp/z$ having root systems $E_8E_6$ and $D_{13}$. Fix $\rho$ to
be a primitive norm 0 vector corresponding to a lattice with root
system $D_{13}$.  As $D_{13}$ has rank one less than the corresponding
lattice, there is a group $\Z$ of automorphisms of the fundamental
domain fixing $z$.  This group $\Z$ has finite index in the full
automorphism group of the fundamental domain, so the full automorphism
group of the fundamental domain must fix $z$. However $z$ is not quite
a Weyl vector, as it has zero inner product with some simple roots
(forming an affine $D_{13}$ Dynkin diagram) and has  inner
product 1 with the others.

There are 10 lattices in the genus $II_{0,12}(3^{+6})$ [S-V].
One is the Coxeter-Todd lattice with no roots, and the others all have
root systems of rank 12.

There are also some automorphic forms of singular weight corresponding
to the following lattices and reflective forms:
\item{1} 
$II_{2,14}(3^{-8})$, $\eta_{1^{-6}3^{-6}}$.
\item{2} 
$II_{2,8}(3^{+7})$, $\eta_{1^{3}3^{-9}}$.
\item{3} 
$II_{2,8}(3^{+3})$, $3^2\eta_{1^{-9}3^{3}}$.

The lattice $II_{2,8}(3^{+5})$ appears in [A-C-T], where it is the
underlying integral lattice of a unimodular Eisenstein lattice. The
automorphic forms for this case have been studied in great detail by
Freitag in [A-F], [F99]. In particular there is one of weight 12
(coming from the function $27\eta_{1^{-9}3^{3}}$) whose restriction to
complex hyperbolic space $\C H^4$ vanishes (to order 3) exactly along
the reflection hyperplanes of a certain complex reflection group
related to the moduli space of cubic surfaces. So its cube root is an
automorphic form of weight 4 with order 1 zeros along all the
reflection hyperplanes [A].

We have seen above that sometimes the Weyl vector of a reflective form
unexpected vanishes because all the singularities just happen to
cancel out. Another way that the Weyl vector can unexpectedly vanish
is if the vectors corresponding to the singularities happen no to
exist (usually when the $p$-rank of $A$ is small).  For example, for 
the lattice
$L=II_{2,8}(3^{-1})$, the automorphic form is constant even though the
vector valued modular form has non-trivial singularities.  The
singularities of the vector valued modular form imply that automorphic
form has zeros corresponding to all norm $4/3$ vectors of $L'$, but
$L'$ happens to have no such vectors so the automorphic form is
constant.

\proclaim $N=4$.

The group $\Gamma_0(4)$ has 3 cusps which can be taken as
$i\infty$ (of width 1) and 0 (of width 4) and $1/2$ (of width 1).  It
has no elliptic points and has genus 0.

\vbox{\halign{
#~&~#~&~#~&~#~&~#~&~#~\cr
Group         &index&$\nu_2$&$\nu_3$&$\nu_\infty$&genus\cr
$\Gamma_0( 4)$&    6&      0&      0&           3&    0\cr
}}
\halign{
~#~&~#~&~#~&~$#$~&~#~&~#~&~#~&~#~\cr
cusps     & width & characters&\eta               & zero & weight &character\cr
      1=0 &    4  &&1^{ 8}2^{-4}              &  1   & 2      &\cr
      1/2 &    1 &$\chi_\theta=-i$ &1^{-2}2^{ 5}4^{-2}&1/4&1/2&$\chi_\theta$\cr
$1/4=i\infty$ &1  &&      2^{-4}4^{ 8}        &  1   & 2     & \cr
}

The double cover of
$\Gamma_0(4)$ is the product of its center of order 4 (generated by
$Z$) and a subgroup that can be identified with $\Gamma_1(4)$.

The ring of modular forms of integral or half integral weight
for $\Gamma_1(4)$ is a polynomial ring generated by
$\theta_{A_1}(\tau) = 1+2q+2q^4+O(q^9)$ of weight $1/2$ and
$\eta(\tau)^{8}\eta(2\tau)^{-4}= 1-8q+\cdots$
of weight 2.
The Hilbert function is
$1/(1-u_\theta x^{1/2})(1-x^2)$. 
The ideals of cusp forms vanishing at $i\infty$, 0, or $1/2$ are generated
by $\eta(2\tau)^{-4}\eta(4\tau)^{8}$, $\eta(\tau)^{8}\eta(2\tau)^{-4}$,
$\eta(\tau)^{-2}\eta(2\tau)^{5}\eta(4\tau)^{-2}$. Note that the last function
has a zero of order $1/4$ at $1/2$.
The ideal of cusp forms of even weight is generated by
$\Delta_{4+}(\tau)=\eta(2\tau)^{12}$
 of weight 6.

If $L$ is a unimodular positive definite lattice
then $\theta_L(2\tau)$ is a modular form for $\Gamma_1(4)$.

Next we find some reflective singularities. To reduce the number of
cases to consider we will
assume that $A=L'/L$ has exponent 2, so that $A$ is
$II(2_t^{+n})$ for some $t$ and $n$.
As usual, poles of order 1 are reflective singularities at the cusps
$i\infty$ and 0. At the cusp $1/2$, poles of order $1/4$ and $1/2$
are reflective, because all elements of $A^{2*}$ have order 1 or 2.
If the parity vector of $A$ does not have norm 0 mod 2 then
a pole of order $1$ at $1/2$ is also reflective.
Finally, at the cusp $0$ poles of order 1 or 2 are reflective, and if
$A$ has no non-zero vectors of norm $0\bmod 2$ then poles of order 4
are reflective.

The form $\eta(\tau)^{-12}\eta(2\tau)^{-2}\eta(4\tau)^{4}$ shows that
all level 4  exponent 2 lattices of signature $-14$ have non-zero reflective
modular forms. By multiplying this form by $\theta_{A_1}(\tau)^n$ for
$n\ge 1$ we see that the level 4 lattices of signature at least $-14$
have non-zero reflexive modular forms.

We can find many examples of eta quotients that are eigenforms of Hecke
operators by finding eta quotients with poles of order at most 1
at all cusps. This gives 15 non-constant examples as follows:
$\eta_{1^{-2} 2^{  5} 4^{-2}}$,
$\eta_{1^{-4} 2^{ 10} 4^{-4}}$,
$\eta_{1^{-8} 2^{ 20} 4^{-8}}$,
$\eta_{1^{ 8} 2^{- 4}       }$,
$\eta_{1^{ 6} 2^{  1} 4^{-2}}$,
$\eta_{1^{ 4} 2^{  6} 4^{-4}}$,
$\eta_{       2^{ 16} 4^{-8}}$,
$\eta_{       2^{- 4} 4^{ 8}}$,
$\eta_{1^{-2} 2^{  1} 4^{ 6}}$,
$\eta_{1^{-4} 2^{  6} 4^{ 4}}$,
$\eta_{1^{-8} 2^{ 16}       }$,
$\eta_{1^{-8} 2^{  8} 4^{-8}}$,
$\eta_{1^{-6} 2^{  3} 4^{-6}}$,
$\eta_{1^{-4} 2^{- 2} 4^{-4}}$,
$\eta_{       2^{-12}       }$.
The inverses of these forms are often reflective forms for various lattices.
Note that the forms $\eta_{1^{-6}2^{15}4^{-6}}$, $\eta_{1^{ 2}2^{11}4^{-6}}$,
$\eta_{1^{-6}2^{11}4^{ 2}}$, $\eta_{1^{ 2}2^{ 7}4^{ 2}}$ are
eigenfunctions of Hecke operators, but as they have a zero of order $3/4$
at $1/2$ their inverses do not usually give reflective automorphic forms
(except for rather special discriminant forms).
The form $\eta_{1^{ 2}2^{ 7}4^{ 2}}$ is the highest weight
eta product I know of that is an eigenform and has non-integral weight.

Most of the time lattices of positive signature with reflective forms
do not seem to be  interesting, but there are some exceptions. For example,
there is a reflective
form for the lattice $II_{2,1}(2^{+1}_1)$. The corresponding automorphic
form is essentially $E_6$, which is the denominator function of a
generalized Kac-Moody algebras. See [B95, section 15, example 2] 
for more details.

The lattices
$II_{1,19}(2^{+2}_{6})$
$II_{1,15}(2^{+2}_{2})$
$II_{1,11}(2^{+2}_{6})$
have reflective forms of type $\theta_{D_n}/\Delta$. They are the even
sublattices of odd unimodular lattices, and have cofinite reflection
groups, as was first found by Vinberg ([V75]).

The function $\theta_{E_7}/\Delta$ is a reflective
form for the lattices $II_{n,n+17}(2^{+1}_{-1})$. In particular we
find the  Nikulin's example of the Lorentzian lattice $II_{1,18}(2^{+1}_{-1})$
whose reflection group is arithmetic. This example can also be
constructed as the orthogonal complement of an $E_7$ in $II_{1,25}$.

Yoshikawa [Y] used the automorphic forms coming from
the modular forms
$\eta(\tau)^{-8}\eta(2\tau)^8\eta(4\tau)^{-8}\theta_{A_1}(\tau)^k$
to construct automorphic products (for odd unimodular lattices).
These automorphic products
are the squares of discriminant forms of various
moduli spaces of ``generalized Enriques surfaces'', and can also
be constructed
using analytic torsion.

\proclaim $N=5$.

\vbox{\halign{
#~&~#~&~#~&~#~&~#~&~#~\cr
Group         &index&$\nu_2$&$\nu_3$&$\nu_\infty$&genus\cr
$\Gamma_0( 5)$&    6&      2&      0&           2&    0\cr
}}
\halign{
~#~&~#~&~$#$~&~#~&~#~&~$#$~\cr
cusps     & width & \eta            & zero & weight \cr
        0 &    5  &  1^{ 5}5^{-1}   &  1   & 2      \cr
$i\infty$ &    1  &  1^{-1}5^{ 5}   &  1   & 2      \cr
}

The ring of modular forms for $\Gamma_0^2(5)$
is not a polynomial ring, but is generated by the 3-dimensional
space of weight 2 forms, which is spanned by
$\eta(\tau)^5\eta(5\tau)^{-1}$, $\eta(\tau)^{-1}\eta(5\tau)^5$
(of nontrivial character)
and $E_2(\tau)-5E_2(5\tau)$ (of trivial character).
The Hilbert function is
$(1+x^2)/(1-u_5x^{2})^2$.

Remark. The ring of all modular forms of integral weight for
$\Gamma_1(5)$ is a polynomial ring generated by the  weight 1
Eisenstein series
$1+(3+i)( q + (1 - i)q^2 + (1 + i)q^3 - iq^4 + q^5+O(q^6))$
and its complex conjugate.
These correspond to the two complex conjugate
order 4 characters of $\Z/4\Z$, and each of them has a simple
zero at one of the elliptic points and no other zeros.
The subring of forms of even weight is the ring of
modular forms for
$\Gamma^2_0(5)$.

Even lattices of level 5 all have signature divisible by 4.  The form
$\Delta_5(\tau)^{-1}$ shows that all level 5 even lattices of signature $-8$
and even 5-rank have non-zero reflective modular forms. (So does the lattice
of 5-rank 1; see below.)
If we multiply
this form by products of powers of $E_2(\tau)-5E_2(5\tau)$ and
$\eta(\tau)^5\eta(5\tau)^{-1}$ we also see that all even level 5
lattices of signature at least $-4$ have non-zero reflective modular forms.

For the discriminant forms $A=II(5^{\pm1})$ or $II(5^{+2})$
the only norm 0 element is 0, so $q_5^{-1}$
and $q_5^{-5}$ are all reflective singularities at the
cusp 0.  If we take $A$ to be $II(5^{-1})$ and take the signature to
be $-8$ then there is a reflective form.  This gives a
Lorentzian lattice $II_{1,9}(5^{-1})$ with a reflection group of
finite index.  If we take $L$ to be $II_{1,17}(5^{-1})$ then there is
a reflective automorphic form. (This is slightly surprising as the
space of forms with a pole of order at most 1 at $i\infty$ and a pole of
order at most 5 at 0 is 2 dimensional, so we would normally expect
there to be no such forms satisfying the 2 conditions that the coefficients
of $q_5^{-2}$ and $q_5^{-3}$ both vanish. However it turns out that these
two conditions are not independent; in fact the modular form we
get has ``complex multiplication'', (see [Ri]) meaning that the coefficient
of $q_5^n$ is 0 whenever $n\equiv 2,3\bmod 5$.)
In spite of the existence of a non-zero reflective modular form,
the reflection group of $II_{1,17}(5^{-1})$ is not cofinite,
and does not even have virtually free abelian index.
(In particular, this lattice is a counterexample to several
otherwise plausible conjectures about Lorentzian lattices with non-zero
reflective modular forms.)
As a substitute for this,
the lattice is very closely related to Bugaenko's largest example
of a co-compact hyperbolic reflection group. In fact $II_{1,17}(5^{-1})$
can be made into a lattice over $\Z[\phi]$, and Bugaenko [B] showed that
the corresponding hyperbolic reflection group was co-compact.
The relationship between Bugaenko's reflection group and the reflective
form is  rather mysterious. The lattice has 5 orbits of primitive
norm 0 vectors, corresponding to the 5 elements of the genus
$II_{0,16}(5^{-1})$, which have root systems $A_2A_{14}$, $E_7A_9$,
$E_6D_9$, $E_8E_75A_1$, $D_{14}A_15A_1$.

It is possible to produce some examples of co-compact hyperbolic reflection
group from level 5 lattices as follows.

\proclaim Lemma 12.1. Suppose that $L$ is an even  Lorentzian lattice of
level 5, and suppose that there is a self adjoint endomorphism $\phi$
of $L$ such that $\phi^2=\phi+1$. Let $H^\phi$ be the hyperbolic space of
the Lorentzian eigenspace $(L\otimes \R)^\phi$. Then the subgroup of the
reflection of $L$ acting on $H^\phi$ is a hyperbolic reflection group of
$H^\phi$. If $W$ is cofinite then $W^\phi$ is co-compact.

Proof. Let $H$ be the hyperbolic space of $L$, and $H^\phi$ the
subspace of it fixed by $\phi$. The main point is that the
intersection of any reflection hyperplane of $W$ with $H^\phi$ is a
reflection hyperplane of the group $W^\phi$ acting on $H^\phi$. To see
this, recall that a reflection of $W$ is the reflection of a norm $-2$
vector of $L$ or a norm $-2/5$ vector of $L'$. First suppose first
that $v$ is a norm $-2$ vector of $L$. As $v^\perp$ intersects
$H^\phi$, $v^{\phi\perp} $ and $v$ must generate a negative definite
space.  This easily implies that $v$ and $v^\perp$ span a lattice
isomorphic to $A_1^2$, and the product of two reflections of this
lattice is the automorphism $-1$ which commutes with $\phi$. This is a
reflection of $W^\phi$ acting on $H^\phi$ whose reflection hyperplane
is $v^\perp\cap H^\phi$. The argument when $v$ is a norm $-2/5$ vector of
$L'$ is similar.

It now follows that $W^\phi$ is a reflection group acting on $H^\phi$
whose fundamental domain is the intersection of
$H^\phi$ with a fundamental domain of $W$ acting on $H$.

Finally if $W$ is cofinite then all norm 0 vectors in the fundamental
domain of $W$ are rational and therefore cannot be fixed by $\phi$, so
the fundamental domain of $W^\phi$ has no norm 0 vectors in it and is
therefore compact. This proves lemma 12.1.

Unfortunately, this lemma does not give the largest examples found by Bugaenko.

If we take $A$ to be $II(5^{+1})$ and take the signature to be $-12$
then $q_5^{-4}$ is a reflective singularity at $0$ as $A$
has no nonzero elements of norm $-4/5 \bmod 2$, and
$q_5^{-5}$ is reflective as any norm 0 element of $A$ is 0.
So $A$ has a reflective modular form of weight $-6$, level 5, and
character $\chi_5$ whose singularity at $i\infty$ is a multiple of $q^{-1}$ and
whose singularity at 0 is a linear combination of $q_5^{-1}$ and
$q_5^{-4}$ and $q_5^{-5}$.  (There is a 2 dimensional space of such forms.) 
This gives a Lorentzian
lattice $L=II_{1,13}(5^{+1})$ with $\Aut^+(L)/R$ infinite
dihedral. It is the orthogonal complement of an $A_4$
in $II_{1,17}$. This case is similar to $II_{1,15}(3^{-1})$.

The (level 1) form
$E_6/\Delta$ has a singularity at 0 of the form $q^{-1}=q_5^{-5}$
so it is a reflective modular form for the lattice $L=II_{1,13}(5^{-2})$.
The corresponding vector valued modular form is 0. 
The lattice $L$ is a module over $\Z[\phi]$. It may be the lattice
of the orthogonal complement of an $I_2(5)$ in Bugaenko's lattice,
which would imply that it has a co-compact reflection group.

As in the case $N=3$ we also get a few examples of automorphic forms
of singular weight coming from the reflective forms
$\eta_{1^{-4}5^{-4}}$, $\eta_{1^{-1}5^{ 5}}$, and $\eta_{1^{
5}5^{-1}}$.  There are 5 lattices in the genus $II_{0,8}(5^{+4})$
corresponding to the case $II_{2,10}$ by [S-V].  One has no roots and
by [S-H p. 744] the roots systems of the other 4 are $A_1^45A_1^4$,
$A_2^25A_2^2$, $A_45A_4$, $D_45D_4$.

Problem: does the lattice $II_{2,6}(5^{+3})$ 
correspond to some nice moduli space,
in  the same way that the corresponding lattices 
$II_{2,10}(2^{+2})$ and $II_{2,8}(3^{+5})$ for levels 2 and
3 correspond to the moduli spaces of Enriques surfaces or cubic surfaces?

\proclaim $N=6$.

\vbox{\halign{
#~&~#~&~#~&~#~&~#~&~#~\cr
Group         &index&$\nu_2$&$\nu_3$&$\nu_\infty$&genus\cr
$\Gamma_0( 6)$&   12&      0&      0&           4&    0\cr
}}
\halign{
~#~&~#~&~$#$~&~#~&~#~&~$#$~\cr
cusps     & width & \eta                     & zero & weight \cr
      1=0 &    6  &1^{ 6}2^{-3}3^{-2}6^{ 1}  &  1   & 1      \cr
      1/2 &    3  &1^{-3}2^{ 6}3^{ 1}6^{-2}  &  1   & 1      \cr
      1/3 &    2  &1^{-2}2^{ 1}3^{ 6}6^{-3}  &  1   & 1      \cr
$1/6=i\infty$ &1  &1^{ 1}2^{-2}3^{-3}6^{ 6}  &  1   & 1      \cr
}

The group
$\Gamma_0(6)$ is the product of $\Gamma_1(6)$ and its center of
order 2 generated by $Z$.
The forms of trivial character have even weight, and those of nontrivial
character have odd weight.

The ring of modular forms of integral   weight
for $\Gamma_0(6)$ is a polynomial ring generated by
$E_1(\tau,\chi_3)$ and $E_1(2\tau,\chi_3)$.
The Hilbert function is $1/(1-u_3x)^2$. 
The ideal of cusp forms is generated by
$\Delta_{6+}(\tau)=\eta(\tau)^{2}\eta(2\tau)^2\eta(3\tau)^{2}\eta(6\tau)^2$
 of weight 4.

We can find eta quotients with an given integral order poles
at the cusps. In particular we find the following 15 non-constant holomorphic
eta quotients with
zeros of order at most 1 at all cusps:
$\eta_{1^{ 6} 2^{-3} 3^{-2} 6^{ 1}}$,
$\eta_{1^{-3} 2^{ 6} 3^{ 1} 6^{-2}}$,
$\eta_{1^{-2} 2^{ 1} 3^{ 6} 6^{-3}}$,
$\eta_{1^{ 1} 2^{-2} 3^{-3} 6^{ 6}}$,
$\eta_{1^{ 3} 2^{ 3} 3^{-1} 6^{-1}}$,
$\eta_{1^{-1} 2^{-1} 3^{ 3} 6^{ 3}}$,
$\eta_{1^{ 4} 2^{-2} 3^{ 4} 6^{-2}}$,
$\eta_{1^{-2} 2^{ 4} 3^{-2} 6^{ 4}}$,
$\eta_{1^{ 7} 2^{-5} 3^{-5} 6^{ 7}}$,
$\eta_{1^{-5} 2^{ 7} 3^{ 7} 6^{-5}}$,
$\eta_{1^{ 1} 2^{ 4} 3^{ 5} 6^{-4}}$,
$\eta_{1^{ 4} 2^{ 1} 3^{-4} 6^{ 5}}$,
$\eta_{1^{ 5} 2^{-4} 3^{ 1} 6^{ 4}}$,
$\eta_{1^{-4} 2^{ 5} 3^{ 4} 6^{ 1}}$,
$\eta_{1^{ 2} 2^{ 2} 3^{ 2} 6^{ 2}}$.
Their inverses give numerous examples of reflective forms for various lattices.
For example, $\eta_{1^{-2} 2^{-2} 3^{-2} 6^{-2}}$ is a reflective form
for all even level 6 lattices of signature $-8$, and by multiplying
by a power of $E_1(\tau,\chi_3)$ we get reflective forms whenever
the signature is at least $-8$. The other eta quotients give many examples
where roots of certain norms are excluded.

We can also find many examples of signature less than $-8$ if
we restrict the 2-rank or 3-rank to be at most 2. One example 
of such a lattice with an arithmetic reflection group is
the orthogonal complement $II_{1,15}(2^{-2}3^{-1})$ 
of a $D_4E_6$ root system in $II_{1,25}$.

\proclaim $N=7$.

\vbox{\halign{
#~&~#~&~#~&~#~&~#~&~#~\cr
Group         &index&$\nu_2$&$\nu_3$&$\nu_\infty$&genus\cr
$\Gamma_0( 7)$&    8&      0&      2&           2&    0\cr
}}
\halign{
~#~&~#~&~$#$~&~#~&~#~&~$#$~\cr
cusps     & width & \eta            & zero & weight \cr
        0 &    7  &  1^{ 7}7^{-1}   &  2   & 3      \cr
$i\infty$ &    1  &  1^{-1}7^{ 7}   &  2   & 3      \cr
}

The ring of modular forms of integral weight for $\Gamma^2_0(7)$ is
generated by $E_1(\tau,\chi_7)$ (of weight 1 which vanishes at both
elliptic points), $\Delta_{7+}(\tau)$ (of weight 3 which vanishes
at both cusps), and the two weight 3 Eisenstein series.
The Hilbert function is
$(1+u_7x^3)/(1-u_7x)(1-u_7x^3)$. 
The ideal of cusp forms is generated by
$\Delta_{7+}(\tau)=\eta(\tau)^{3}\eta(7\tau)^3$ of weight 3.  Note
that the ideal of forms vanishing at $i\infty$ is not principal.  The
function $\eta(\tau)^4\eta(7\tau)^{-4}$ is a Hauptmodul for
$\Gamma_0(7)$.

Remark. The ring of modular forms for $\Gamma_1(7)$ of integral weight
has a simpler structure: it is generated by the three weight 1 forms
$E_1(\tau,\chi_7)$, $E_1(\tau,\chi)$, $E_1(\tau,\bar\chi)$, where
$\chi$ is a character of $\Z/7\Z$ of order 6 and $\bar\chi$ is its
complex conjugate. The ideal of relations between these generators is
generated by $E_1(\tau,\chi_7)^2-E_1(\tau,\chi)E_1(\tau,\bar\chi)$.
We can even embed this  into a polynomial ring of modular forms: all
the zeros of the forms $E_1(\tau,\chi)$ and $E_1(\tau,\bar\chi)$ have
order 2, so their square roots are also modular forms (of half
integral weight for a strange character of $\Gamma_0(7)$), and they
generate a polynomial ring whose elements of integral weight are the
modular forms for $\Gamma_1(7)$.

Any even lattice of level 7 and signature at least $-6$
has a reflective form of the form $E_1(\tau,\chi_7)^n\eta_{1^{-3}7^{-3}}$.

The automorphic form associated to $\eta_{1^{-3}7^{-3}}$ and the
lattice $II_{2,8}(7^{+5})$ has singular weight and is the denominator
function of a generalized Kac-Moody algebra. The reflection group
of $II_{1,7}(7^{-3})$ has a norm 0 Weyl vector. The corresponding genus
$II_{0,6}(7^{-3})$ has 3 elements [S-H, proposition 3.4a], and by [S-H
table 1] there is one with no roots (corresponding to the norm 0 Weyl
vector), one with root system $A_37A_3$ and one with root system
$A_1^37A_1^3$.

For the discriminant form $II(7^{-1})$ the singularities $q_7^{-1}$
and $q_7^{-7}$ are reflective.  The lattice $II_{1,11}(7^{-1})$ has
a reflective form, and is the orthogonal complement of an $A_6$
in $II_{1,17}$. The quotient $\Aut^+(L)/R$ is infinite dihedral, and
fixes a norm 0 vector corresponding to a lattice in the genus
$II_{0,10}(7^{-1})$ with root system $D_9$. There is a second lattice
in this genus, isomorphic to the sum of $E_8$ and a 2 dimensional
definite lattice of determinant 7, so its root system is $E_8A_17A_1$.

\proclaim $N=8$.

\vbox{\halign{
#~&~#~&~#~&~#~&~#~&~#~&~#~\cr
Group         &index&$\nu_2$&$\nu_3$&$\nu_\infty$&genus\cr
$\Gamma_0( 8)$&   12&      0&      0&           4&    0\cr
}}
\halign{
~#~&~#~&~#~&~$#$~&~#~&~#~&~#~\cr
cusps     & width &characters        & \eta                     & zero & weight
&character           \cr
      1=0 &    8  &                  &1^{ 4}2^{-2}              &  1   & 1    
&$\chi_\theta^2$     \cr
      1/2 & 2 &$\chi_2=\chi_\theta=-1$ &1^{-2}2^{ 5}4^{-2}      &  1/2 & 1/2 
&$\chi_\theta$       \cr
      1/4 &    1  &$\chi_2=-1$       &      2^{-2}4^{ 5}8^{-2}  &  1/2 & 1/2  
&$\chi_\theta\chi_2$ \cr
$1/8=i\infty$ &1  &                  &            4^{-2}8^{ 4}  &  1   & 1    
&$\chi_\theta^2$     \cr
}

The double cover of $\Gamma_0(8)$ is the product of its center
of order 4 (generated by $Z$) and a subgroup that can be identified
with its image $\Gamma_0(8)\cap \Gamma_1(4)$.
The ring of modular forms of integral or half integral weight for
$\Gamma^2_0(8) = \Gamma_1(8)$ is a polynomial ring generated by
$\eta(2\tau)^{-2}\eta(4\tau)^{5}\eta(8\tau)^{-2}$ (of weight $1/2$ and
character $\chi_2\chi_\theta$) and
$\eta(\tau)^{-2}\eta(2\tau)^{5}\eta(4\tau)^{-2}$ (of weight $1/2$ and
character $\chi_\theta$).  The Hilbert function is
$1/(1-u_\theta u_{2}x^{1/2})(1-u_\theta x^{1/2})$.

There are many level 8 discriminant forms and many possible reflective
singularities.  Together they give a bewildering number of examples of
level 8 lattices with reflective forms; they are probably best left to
a computer to classify.  As examples we will just mention
$II_{1,18}(4^{+1}_{ 7})$, $II_{1,16}(4^{+1}_{ 1})$,
$II_{1,14}(4^{-1}_{ 3})$, $II_{1,12}(4^{-1}_{ 5})$.  These are the
even sublattices of some of the odd unimodular lattices with cofinite
reflection groups found by Vinberg and Kaplinskaja [V-K].

\proclaim $N=9$.

\vbox{\halign{
#~&~#~&~#~&~#~&~#~&~#~\cr
Group         &index&$\nu_2$&$\nu_3$&$\nu_\infty$&genus\cr
$\Gamma_0( 9)$&   12&      0&      0&           4&    0\cr
}}
\halign{
~#~&~#~&~$#$~&~#~&~#~&~$#$~\cr
cusps     & width & \eta                     & zero & weight \cr
      1=0 &    9  &1^{ 3}3^{-1}              &  1   & 1      \cr
 1/3, 2/3 &    1  &1^{-3}3^{10}9^{-3}        &  1,1 & 2      \cr
$1/9=i\infty$ &1  &      3^{-1}9^{ 3}        &  1   & 1      \cr
}

The group $\Gamma_0(9)$ is the product of its center of order 2
and the group $\Gamma_0^2(9)$, so a form has nontrivial character if and only
if it has odd weight.
The ring of modular forms of integral  weight
for $\Gamma_0^2(9)$ is a polynomial ring  generated by
$\eta(\tau)^{3}\eta(3\tau)^{-1}$ and
$\eta(9\tau)^{3}\eta(3\tau)^{-1}$.
The Hilbert function is $1/(1-x)^2$.

For the sake of completeness we also describe generators for
the ring of all integral weight
modular forms for $\Gamma_1(9)$. This is a 3 dimensional free module over
the ring of modular forms for $\Gamma^2_0(9)$, with a basis consisting
of 1 and the two weight 1 Eisenstein series for the two order 6
characters of $\Z/6\Z$. Each of these Eisenstein series
has zeros of order $1/3$ and $2/3$ at the cusps $1/3$ and $2/3$
(not necessarily in that order). Note that the character of
a modular form can be read off from the parity of its weight and
the fractional part of the order of the zero at $1/3$.
We can embed this ring in a polynomial ring, generated by the cube roots
of the two weight 1 modular forms with poles of order 1 at $1/3$ or $2/3$.

The form $\eta(\tau)^{-3}\eta(3\tau)^2\eta(9\tau)^{-3}$
is a reflective form for all even level 9 lattices of signature
$-4$, and by multiplying it by a suitable power of (say) $\theta_{A_2}(\tau)$
we get reflective forms whenever the signature is at least $-4$.

A few examples of modular forms that might correspond to automorphic forms
of singular weight on some lattices are
$\eta_{       3^{-8}       }$,
$\eta_{1^{-3} 3^{ 1}       }$,
$\eta_{       3^{ 1} 9^{-3}}$, and
$\eta_{1^{-3} 3^{ 2} 9^{-3}}$.

\proclaim $N=10$.

\vbox{\halign{
#~&~#~&~#~&~#~&~#~&~#~\cr
Group         &index&$\nu_2$&$\nu_3$&$\nu_\infty$&genus\cr
$\Gamma_0(10)$&   18&      2&      0&           4&    0\cr
}}
\halign{
~#~&~#~&~$#$~&~#~&~#~&~$#$~\cr
cusps     & width & \eta                      & zero & weight \cr
      1=0 &   10  &1^{10}2^{-5}5^{-2}10^{ 1}  &  3   & 2      \cr
      1/2 &    5  &1^{-5}2^{10}5^{ 1}10^{-2}  &  3   & 2      \cr
      1/5 &    2  &1^{-2}2^{ 1}5^{10}10^{-5}  &  3   & 2      \cr
$1/10=i\infty$&1  &1^{ 1}2^{-2}5^{-5}10^{10}  &  3   & 2      \cr
}

The ring of modular forms of integral   weight
for $\Gamma^2_0(10)$ is  generated by the 7 dimensional space of forms
of weight 2, and the ideal of relations between them is
generated by 15 quadratic relations.
The Hilbert function is $(1+5x^2)/(1-x)^2$.

Remark. We can embed the ring of modular forms for $\Gamma^2_0(10)$ in a
polynomial ring as follows.  The ring of modular forms of integral
weight for $\Gamma_1(10)$ is generated by the 4 dimensional space of
forms of weight 1. We can find two of these forms that have zeros of
order 3 at the two elliptic points of order 2. Their cube roots are
modular forms of weight $1/3$ for characters of order 3 of
$\Gamma_1(10)$, and generate a polynomial ring in 2 variables.  The
space of modular forms for $\Gamma^2_0(10)$ can be identified with the
polynomials of degree divisible by 6.

The cube root of the product of any three of the forms above
with zeros only at one cusp is a form with a zero
of order 1 at 3 of the 4 cusps, and these 4 forms are a basis
of the space of weight 2 forms with non-trivial character.
Their inverses are the functions
$\eta_{1^{-1} 2^{-2} 5^{-3} 10^{ 2}}$,
$\eta_{1^{-2} 2^{-1} 5^{ 2} 10^{-3}}$,
$\eta_{1^{-3} 2^{ 2} 5^{-1} 10^{-2}}$,
$\eta_{1^{ 2} 2^{-3} 5^{-2} 10^{-1}}$.

Any one of these 4 functions (of non-trivial character and weight $-2$)
shows that any even
level 10 lattice of signature $-4$ and odd 5-rank has a reflective
form. There is also a weight $-2$ form of trivial character whose
poles and zeros are a pole of order 1 at each cusp and order 1 zeros
at the two elliptic points.
(Construction: take a linear combination
of the weight 2 Eisenstein series with trivial character that vanishes
at 2 cusps (this automatically vanishes at the two elliptic points),
then divide it by an eta product with order 2 zeros at these cusps
and order 1 zeros at the other two cusps.) 
This is a reflective form
for the even level 10 lattice of signature $-4$ with even 5-rank.

So every even level 10 lattice of signature $\ge -4$ has a reflective form.

\proclaim $N=11$.

\vbox{\halign{
#~&~#~&~#~&~#~&~#~&~#~\cr
Group         &index&$\nu_2$&$\nu_3$&$\nu_\infty$&genus\cr
$\Gamma_0(11)$&   12&      0&      0&           2&    1\cr
}}
\halign{
~#~&~#~&~$#$~&~#~&~#~&~$#$~\cr
cusps     & width & \eta            & zero & weight \cr
        0 &   11  &1^{ 11}11^{-1}   &  5   & 5      \cr
$i\infty$ &    1  &1^{-1}11^{ 11}   &  5   & 5      \cr
}

The ring of modular forms of integral weight for $\Gamma^2_0(11)$ is
generated by $E_1(\tau,\chi_{11})$ (of weight 1 which vanishes at both
elliptic points), $\Delta_{11+}(\tau)$ (of weight 2 which vanishes
at both cusps), and the two weight 3 Eisenstein series.
The Hilbert function is
$(1+x^3)/(1-x)(1-x^2)$. 
The ideal of cusp forms is generated by
$\Delta_{11+}(\tau)=\eta(\tau)^{2}\eta(11\tau)^2$ of weight 2. The
ideal of forms vanishing at $i\infty$ is not principal.
The function  $\eta(\tau)^{12}\eta(11\tau)^{-12}$
has a pole of order 5 at $i\infty$ and a zero of order 5 at 0
and is a modular function for $\Gamma_0(11)$, showing that
$0$ is a torsion point of order 5 on the modular elliptic curve of
$\Gamma_0(11)$. (In fact this point generates the subgroup of rational
points on this elliptic curve.)

The forms $E_1(\tau,\chi_{11})^n\eta_{1^{-2} 11^{-2}}$ show that
all even lattices of level 11 and signature at least $-4$ are reflective.
The lattice $II_{1,7}(11^{-1})$ has a reflection group of infinite
dihedral index in its automorphism  group. The corresponding
genus $II_{0,6}(11^{-1})$
contains just one lattice, which has root system $D_5$.

\proclaim $N=12$.

\vbox{\halign{
#~&~#~&~#~&~#~&~#~&~#~\cr
Group         &index&$\nu_2$&$\nu_3$&$\nu_\infty$&genus\cr
$\Gamma_0(12)$&   24&      0&      0&           6&    0\cr
}}
\halign{
~#~&~#~&~#~&~$#$~&~#~&~$#$~\cr
cusps     & width &characters& \eta                     & zero & weight \cr
      1=0 &   12  && 1^{ 6}2^{-3}3^{-2}6^{ 1} &  2   & 1      \cr
      1/2 &    3  &$\chi_\theta=i$& 1^{-6}2^{15}3^{2}4^{-6}6^{-5}12^{2}& 2 
&1\cr
      1/3 &    4  && 1^{-2}2^{ 1}3^{ 6}6^{-3} &  2   & 1      \cr
      1/4 &    3  && 2^{-3}4^{ 6}6^{ 1}12^{-2}&  2   & 1      \cr
      1/6 &    1  &$\chi_\theta=-i$& 1^{2}2^{-5}3^{-6}4^{2}6^{15}12^{-6}& 2 
&1\cr
$1/12=i\infty$&1  && 2^{1 }4^{-2}6^{-3}12^{ 6}&  2   & 1      \cr
}

The ring of modular forms is generated by the forms
$\theta(\tau)$, $\theta(3\tau)$, $E_1(\tau,\chi_3)$, and $E_1(2\tau,\chi_3)$
and the Hilbert function is $(1+x^{1/2}+2x)/(1-x^{1/2})(1-x)$.

There are quite a lot of modular forms whose zeros are all zeros of
order at most 1 at cusps: we can find forms with zeros of order $1/2$
at $1/2$ and $1/6$ and an odd number of zeros at the other cusps, or
forms whose zeros at $1/2$ and $1/6$ have orders $(0,0)$, $(1/4,3/4)$,
$(3/4,1/4)$, or $(1,1)$ and that have an even number of zeros at the
other cusps. The maximum weight of these forms is 3, attained by the
form $\eta(2\tau)^3\eta(6\tau)^3$ with a zero of order 1 at every
cusp.
The inverses of these forms are reflective forms for many lattices of 
signature up to $-6$.

\proclaim $N=13$.

\vbox{\halign{
#~&~#~&~#~&~#~&~#~&~#~\cr
Group         &index&$\nu_2$&$\nu_3$&$\nu_\infty$&genus\cr
$\Gamma_0(13)$&   14&      2&      2&           2&    0\cr
}}
\halign{
~#~&~#~&~$#$~&~#~&~#~&~$#$~\cr
cusps     & width & \eta            & zero & weight \cr
        0 &   13  &1^{ 13}13^{-1}   &  7   & 6      \cr
$i\infty$ &    1  &1^{-1}13^{ 13}   &  7   & 6      \cr
}

The Hilbert function is $(1+2x^2+6x^4+5x^6)/(1-x^2)(1-x^6)$.
The space of weight 2 forms for $\Gamma^2_0(13)$ is 3 dimensional,
spanned by $E_2(\tau)-13E_3(13\tau)$, $E_2(\tau, \chi_{13})$,
and the cusp form
$E_1(\tau,\chi)^2-E_2(\tau,\bar\chi)^2$ where $\chi$ is an order
4 character of $(\Z/13\Z)^*$.
The ring of modular forms for $\Gamma_1(13)$ is generated by the 6
dimensional space of forms of weight 1, which has a basis of the 6
forms $E_1(\tau,\chi)$ as $\chi$ runs through the 6 odd characters of
$(\Z/13\Z)^*$. Each of these weight 1 Eisenstein series has a zero at an
elliptic point of order 2 and 2 zeros at elliptic points of order 3
and no other zeros.

The function $\eta_{1^2 13^{-2}}$ is a Hauptmodul for $\Gamma_0(13)$.
There is also a Hauptmodul for $\Gamma_0(13)+$. These give automorphic forms
for the lattice $II_{2,2}(13^{+2})$, which are the denominator functions
for generalized Kac-Moody algebras related to elements of order
13 in the monster group.

There are no modular forms of negative weight with poles of order at
most 1 at the cusps, as can be seen from the relation (number of
zeros) = weight$\times$index/12 = weight$\times$7/6 and the fact that
the weight is even and there are only 2 cusps.  The space of cusp
forms of weight 4 and character $\chi_{13}$ has dimension 2, and as
this is the space of obstructions to finding a form of weight $-2$ and
character $\chi_{13}$ with given singularities, we see that there is a
nonzero form of weight $-2$ and character $\chi_{13}$ whose
singularities are a pole of order 1 at
$i\infty$ and a singularity at 0 with terms involving only
$q^{-1}$ and $q^{-3}$. This is a reflective form for the lattices
$II_{n,4+n}(13^{+1})$.

The lattice $L=II_{1,5}(13^{+1})$ is one of the lattices with
$\Aut(L)/R(L)$ infinite dihedral. This is
the orthogonal complement of an $A_{12}$ in
$II_{1,17}$.  There is a unique lattice in the genus
$II_{0,4}(13^{+1})$, and it has root system is $D_3$. This can be seen
from the fact that all such lattices are the orthogonal complement of
an $A_{12}$ in an even 16 dimensional self dual negative definite
lattice.

\proclaim $N=14$.

 The group $\Gamma_0(14)$ is the product of $\Gamma_0^2(14)$ and its
center of order 2 generated by $Z$.

\vbox{\halign{
#~&~#~&~#~&~#~&~#~&~#~\cr
Group         &index&$\nu_2$&$\nu_3$&$\nu_\infty$&genus\cr
$\Gamma_0(14)$&   24&      0&      0&           4&    1\cr
}}
\halign{
~#~&~#~&~$#$~&~#~&~#~&~$#$~\cr
cusps     & width & \eta                      & zero & weight \cr
      1=0 &   14  &1^{14}2^{-7}7^{-2}14^{ 1}  &  6   & 3      \cr
      1/2 &    7  &1^{-7}2^{14}7^{ 1}14^{-2}  &  6   & 3      \cr
      1/7 &    2  &1^{-2}2^{ 1}7^{14}14^{-7}  &  6   & 3      \cr
$1/14=i\infty$&1  &1^{ 1}2^{-2}7^{-7}14^{14}  &  6   & 3      \cr
}

The ring of modular forms of integral   weight
for $\Gamma_0^2(14)$ is  generated by
$E_1(\tau,\chi_7)$, $E_1(2\tau, \chi_7)$, and
$\Delta_{14+}(\tau)$.
The Hilbert function is $(1+x^2)/(1-x)^2$
The ideal of cusp forms is generated by
$\Delta_{14+}(\tau)=\eta(\tau)\eta(2\tau)\eta(7\tau)\eta(14\tau)$
of weight 2.
We can construct some automorphic forms on lattices of signature $-2$ or $-4$
of singular weight
from the modular forms 
$\eta_{1^{-2} 2^{ 1} 7^{-2} 14^{ 1}}$,
$\eta_{1^{ 1} 2^{-2} 7^{ 1} 14^{-2}}$, and
$\eta_{1^{-1} 2^{-1} 7^{-1} 14^{-1}}$.

\proclaim $N=15$.

\vbox{\halign{
#~&~#~&~#~&~#~&~#~&~#~\cr
Group         &index&$\nu_2$&$\nu_3$&$\nu_\infty$&genus\cr
$\Gamma_0(15)$&   24&      0&      0&           4&    1\cr
}}
\halign{
~#~&~#~&~$#$~&~#~&~#~&~$#$~\cr
cusps     & width & \eta                      & zero & weight \cr
      1=0 &   15  &1^{15}3^{-5}5^{-3}15^{ 1}  &  8   & 4      \cr
      1/3 &    5  &1^{-5}3^{15}5^{ 1}15^{-3}  &  8   & 4      \cr
      1/5 &    3  &1^{-3}3^{ 1}5^{15}15^{-5}  &  8   & 4      \cr
$1/15=i\infty$&1  &1^{ 1}3^{-3}5^{-5}15^{15}  &  8   & 4      \cr
}

This case seems very similar to the case $N=14$.
The ring of modular forms of integral weight
for $\Gamma^2_0(15)$ is  generated by
$E_1(\tau,\chi_3)$, $E_1(5\tau,\chi_3)$, and $\Delta_{15+}(\tau)$.
The Hilbert function is $(1+x^2)/(1-x)^2$.

The ideal of cusp forms is generated by
$\Delta_{15+}(\tau)=\eta(\tau)\eta(5\tau)\eta(3\tau)\eta(15\tau)$
 of weight 2.
We can construct some automorphic forms on lattices of signature $-2$ or $-4$
of singular weight
from the modular forms 
$\eta_{1^{-2} 3^{ 1} 5^{ 1} 15^{-2}}$,
$\eta_{1^{ 1} 3^{-2} 5^{-2} 15^{ 1}}$, and
$\eta_{1^{-1} 3^{-1} 5^{-1} 15^{-1}}$.

\proclaim $N=16$.

\vbox{\halign{
#~&~#~&~#~&~#~&~#~&~#~&~#~&~#~\cr
Group         &index&$\nu_2$&$\nu_3$&$\nu_\infty$&genus\cr
$\Gamma_0(16)$&   24&      0&      0&           6&    0\cr
}}
\halign{
~#~&~#~&~#~&~$#$~&~#~&~#~&~#~\cr
cusps     & width &characters & \eta             & zero & weight&character \cr
      1=0 &   16  &           &1^{ 2}2^{-1}      &  1   & 1/2 &$\chi_\theta$\cr
      1/2 &    4  &           &1^{-2}2^{ 5}4^{-2}&  1   & 1/2 &$\chi_\theta$\cr
1/4, 3/4&1&$\chi_2=-1$&2^{-2}4^{ 5}8^{-2}&1/2,1/2&1/2 &$\chi_2\chi_\theta$ \cr
      1/8 &    1  &           &4^{-2}8^{ 5}16^{-2}& 1   & 1/2 &$\chi_\theta$\cr
$1/16=i\infty$&1  &           &     8^{-1}16^{ 2}&  1   & 1/2 &$\chi_\theta$\cr
}

Note that the isomorphism class of a cusp $a/c$ is no longer always determined
by $(c,16)$.

The ring of modular forms of integral or half integral weight
for $\Gamma^2_0(16)$ is  generated by
the 3-dimensional space of forms of weight $1/2$,
which is spanned by  the 5 forms listed above.
The Hilbert function is
$(1+\chi_2\chi_\theta x^{1/2})/(1-\chi_\theta x^{1/2})^2$.
As in the case $N=8$ there seem to be rather a lot of examples.

\proclaim $N=17$.

\vbox{\halign{
#~&~#~&~#~&~#~&~#~&~#~\cr
Group         &index&$\nu_2$&$\nu_3$&$\nu_\infty$&genus\cr
$\Gamma_0(17)$&   18&      2&      0&           2&    1\cr
}}

All modular forms for
$\Gamma_0^2(17)$ have even weight.  The Hilbert function is
$(1+(1+2u_{17})x^2+(3+4u_{17})x^4+x^6)/(1-x^2)(1-x^4)$. 
The group $\Gamma_0(17)+$ has
genus 0, and its Hauptmodul $q^{-1}+7q+14q^2+O(q^3)$ has poles of
order 1 at all cusps. 

The Hauptmodul with poles of order 1 at all cusps is a reflective
modular form for lattices $II_{n,n}(17^{\pm 2m})$.  For example, we
get automorphic forms for the lattice $II_{2,2}(17^{+2})$.  The
automorphic form for $II_{2,2}(17^{+2})$ is the denominator function
of a generalized Kac-Moody algebra associated with an element of order
17 of the monster group. We get similar statements if we replace 17 by
any of the other primes $p=2$, 3, 5, 7, 11, 13, 17, 19, 23, 29, 31,
41, 47, 59, or $71$ such that $\Gamma_0(p)+$ has genus 0.

\proclaim $N=18$.

\vbox{\halign{
#~&~#~&~#~&~#~&~#~&~#~\cr
Group         &index&$\nu_2$&$\nu_3$&$\nu_\infty$&genus\cr
$\Gamma_0(18)$&   36&      0&      0&           8&    0\cr
}}
\halign{
~#~&~#~&~$#$~&~#~&~#~&~$#$~\cr
cusps     & width & \eta                             & zero & weight \cr
      1=0 &   18  &1^{ 6}2^{-3}3^{-2}6^{ 1}          &  3   & 1      \cr
      1/2 &    9  &1^{-3}2^{ 6}3^{ 1}6^{-2}          &  3   & 1      \cr
 1/3, 2/3 &    2  &1^{-6}2^{3}3^{20}6^{-10}9^{-6}18^3&  3,3 & 2      \cr
 1/6, 5/6 &    1  &1^{3}2^{-6}3^{-10}6^{20}9^318^{-6}&  3,3 & 2      \cr
$1/9$     &    2  &      3^{-2}6^{1}9^{ 6}18^{-3}    &  3   & 1      \cr
$1/18=i\infty$&1  &      3^{ 1}6^{-2}9^{-3}18^6      &  3   & 1      \cr
}

The Hilbert function is $(1+2u_3x)/(1-u_3x)^2$, and the ring of
modular forms is generated by the 4-dimensional space of forms of
weight 1, which is spanned by $E_1(\tau,\chi_3)$, $E_1(2\tau,\chi_3)$,
$E_1(3\tau,\chi_3)$, and $E_1(6\tau,\chi_3)$.

There are many eta quotients with poles of order 1 at all cusps: for
any set of 3 or 6 cusps (with multiplicities) such that cusps with the
same denominator have the same multiplicity there is an eta quotient
of weight 1 or 2 with these zeros. This gives 12 such eta quotients of
weight 1 and 8 of weight 2. There are 4 of weight 1 with no poles
at the cusps $1/3$, $2/3$, $1/6$, or $5/6$, and the inverses of these
are reflective forms for even level 18 lattices of signature $-2$.

\proclaim $N=20$.

\vbox{\halign{
#~&~#~&~#~&~#~&~#~&~#~\cr
Group         &index&$\nu_2$&$\nu_3$&$\nu_\infty$&genus\cr
$\Gamma_0(20)$&   36&      0&      0&           6&    1\cr
}}
\halign{
~#~&~#~&~#~&~$#$~&~#~&~#~&~$#$~\cr
cusps&width&characters      & \eta                            &zero&weight\cr
 1=0 &  20 &                &1^{ 10}2^{-5}5^{-2}10^{ 1}              &6 &2\cr
 1/2 &   5 &$\chi_\theta=-i$&1^{-10 }2^{25}4^{- 10}5^{2}10^{-5}20^{2}&6 &2\cr
 1/4 &   5 &                &2^{-5}4^{ 10}10^{ 1}20^{-2}             &6 &2\cr
 1/5 &   4 &                &1^{-2}2^{ 1}5^{ 10}10^{-5}              &6 &2\cr
1/10 &   1 &$\chi_\theta=-i$&1^{2}2^{-5}4^{2}5^{- 10}10^{25}20^{- 10}&6 &2\cr
$1/20=i\infty$&1  &         &2^{1 }4^{-2}10^{-5}20^{ 10}             &6 &2\cr
}

The space of weight $1/2$ forms is spanned by $\eta_{1^{-2}2^54^{-2}}$
(zero of order $5/4$ at $1/2$ and order $1/4$ at $1/10$) and
$\eta_{5^{-2}10^520^{-2}}$
(zero of order $1/4$ at $1/2$ and order $5/4$ at $1/10$).
The space of cusp forms of weight 2 and trivial character is
spanned by $\eta_{2^210^2}$; it has zeros of order 1 at all cusps.

The form $1^12^14^{-1}5^{-1}10^120^1$ has zeros of order $1/2$ at
$1/2$ and $1/10$ and zeros of order 1 at $1$ and $1/20$.
The form $1^{-1}2^14^{1}5^{1}10^120^{-1}$ has zeros of order $1/2$ at
$1/2$ and $1/10$ and zeros of order 1 at $1/4$ and $1/5$.
Their inverses are reflective forms for even lattices of
level 20, signature $-2$, and even 5-rank.

\proclaim $N=23$.

\vbox{\halign{
#~&~#~&~#~&~#~&~#~&~#~\cr
Group         &index&$\nu_2$&$\nu_3$&$\nu_\infty$&genus\cr
$\Gamma_0(23)$&   24&      0&      0&           2&    2\cr
}}
\halign{
~#~&~#~&~$#$~&~#~&~#~&~$#$~\cr
cusps     & width & \eta            & zero & weight \cr
        0 &   23  &1^{ 23}23^{-1}   &  22  & 11      \cr
$i\infty$ &    1  &1^{-1}23^{ 23}   &  22  & 11      \cr
}

The ring of modular forms of integral weight for $\Gamma^2_0(23)$ is
generated by $\Delta_{23+}(\tau)=\eta(\tau)\eta(23\tau)$,
and $E_1(\tau, \chi_{23})$,
and one of the two weight 3 Eisenstein series.
The Hilbert function is
$(1+x^3)/(1-x)^2$.

The ideal of cusp forms is generated by
$\Delta_{23+}(\tau)=\eta(\tau)\eta(23\tau)$ of weight 1.
This is the lowest level for which there is a cusp form of weight 1.

Remark. The modular function
$\eta(\tau)^{12}\eta(23\tau)^{-12}$  has a pole of order 11 at
$i\infty$ and a zero of order 11 at 0, and shows that the cusp 0
gives a torsion point of order 11 on the modular abelian surface of
$\Gamma_0(23)$. See [Sh, p. 197] for more about this.

The forms $E_1(\tau,\chi_{23})^n/\Delta_{23+}(\tau)$ show that
all level 23 lattices of signature at least $-2$ have reflective forms.

The automorphic form of the lattice $II_{2,4}(23^{+3})$ and the
function $1/\Delta_{23+}$ has singular weight and is the denominator
function of a generalized Kac-Moody algebra. This generalized
Kac-Moody algebra contains the Feingold-Frenkel rank 3 Kac-Moody
algebra as a subalgebra, and can be used to explain why the root
multiplicities of the Feingold-Frenkel algebra are often given by
values of the partition function. See [Ni] for details.

\proclaim $N=28$.

The eta product 
$1^12^{-1}4^17^114^{-1}28^1$ has weight 1, character $\chi_7$, and its
zeros are
order 1 zeros at the cusps $1/1$, $1/4$, $1/7$, and $1/28$. (It is nonzero
at the cusps $1/2$ and $1/14$.) So its inverse is a reflective form
for any level 28 even lattice of signature $-2$ and odd 7-rank.

\proclaim $N=30$.

Here are some weight 1 eta quotients with order 1 zeros at 6 of the 8 cusps.
$      2^{ 1}3^{ 1}5^{ 1}6^{-1}10^{-1}       30^{ 1}$ 
(nonzero at $1/3$, $1/5$),
$1^{ 1}      3^{-1}5^{-1}6^{ 1}10^{ 1}15^{ 1}       $ 
(nonzero at $1/6$, $1/10$),
$1^{-1}2^{ 1}3^{ 1}5^{ 1}             15^{-1}30^{ 1}$ 
(nonzero at $1/2$, $1/30$),
$1^{ 1}2^{-1}            6^{ 1}10^{ 1}15^{ 1}30^{-1}$
(nonzero at $1/1$, $1/15$).
These are all modular forms for the character $\chi_{15}$.
So any even lattice of level 30 and signature $-2$ and odd 5-rank has
a reflective form. (The 3-rank is automatically odd for any such lattice.)
The maximal order of an automorphism  of the Leech lattice with fixed points
is 30, and 3 of the eta quotients above occur as generalized cycle shapes
of such order 30 automorphisms.

There is no special reason for stopping at $N=30$: there are 
hints that there might be examples for $N$ up to a few hundred. 

\proclaim 13.~Open problems.

We list a few suggestions for further research.

\noindent{\bf Problem 13.1.} 
Find some analogue of reflective forms for other sorts of hyperbolic
reflection group. In particular explain why the complex hyperbolic
reflection groups found by Allcock [A] have underlying integral
lattices with non-zero reflective modular forms.  Is this true of all
complex hyperbolic reflection groups (except perhaps in small
dimensions)?  Is there a relation between the co-compact hyperbolic reflection
groups found by Bugaenko [B] and lattices with reflective forms?

\noindent{\bf Problem 13.2.}
Find some sort of converse theorem that implies that all
``interesting'' lattices of some sort have non-zero reflective forms.
Bruinier [Br] has recently
proved some related converse theorems, showing that certain sorts
of automorphic infinite products always come from modular forms with
singularities.

\noindent{\bf Problem 13.3.}
Are there any other reflection groups of high dimensional lattices
with co-finite volume other than those listed in section 12?

\noindent{\bf Problem 13.4.}
Which of the rank 3 hyperbolic lattices classified by Nikulin [N97]
have reflective forms? Note that in part III of Nikulin's papers
there are many examples where the Weyl vector has negative norm.

\noindent{\bf Problem 13.5.}
Classify all holomorphic eta quotients whose zeros at cusps are all of
order at most 1. This would be suitable for a computer.

\noindent{\bf Problem 13.6.}
Write a computer program to classify the lattices such that the
space of possible reflective singularities is greater than the dimension 
of the space of cusp forms that give obstructions to the existence of
such a singularity.  This would give a large class of examples of
lattices with reflective forms, and the examples in section 12 suggest
that this would include most of them.

\noindent{\bf Problem 13.7.}
The group $\Gamma_0(24^2)$ has 48 cusps, all conjugate under its
normalizer. The form
$1/\eta(24\tau)$ has poles of order 1 at all cusps.
Are there any  lattices for which it is a reflective form?

\proclaim References.

\item{[A]} D. Allcock, The Leech lattice and complex hyperbolic
reflections, 1997 preprint, available from
http://www.math.utah.edu/\hbox{\~{}}allcock
\item{[A-F]} D. Allcock, E. Freitag,
Cubic surfaces and Borcherds products,
1999 preprint, available from 
http://www.rzuser.uni-heidelberg.de/\hbox{\~{}}t91
\item{[A-C-T]}
D. Allcock,  J. A. Carlson,  D. A. Toledo,
A complex hyperbolic structure for moduli of cubic surfaces.
C. R. Acad. Sci. Paris S\'er. I Math. 326 (1998), no. 1, 49--54.
\item{[BBCO]} C. Batut, D. Bernardi, H. Cohen, M. Olivier,
``User's guide to PARI-GP''. This guide and the PARI programs
can be obtained  from 
ftp://megrez.math.u-bordeaux.fr/pub/pari
\item{[B87]}
R. E. Borcherds, Automorphism groups of Lorentzian
lattices. J. Algebra 111 (1987), no. 1, 133--153.
\item{[B90]}
R. E. Borcherds, Lattices like the Leech lattice. J. Algebra 130
(1990), no. 1, 219--234.
\item{[B95]} R. E. Borcherds,  Automorphic forms on
     $O_{s+2,2}(R)$ and infinite products.
     Invent. Math. 120, p. 161-213 (1995)
\item{[B96]} R. E. Borcherds,  The moduli space of Enriques surfaces and the
     fake monster Lie superalgebra. Topology vol. 35 no. 3, 699-710, 1996.
\item{[B98a]} R. E. Borcherds,
Automorphic forms with singularities on Grassmannians,
alg-geom/9609022.  Invent. Math. 132 (1998), no.  3, 491--562.
\item{[B98b]} R. E. Borcherds,
Coxeter groups, Lorentzian lattices, and K3 surfaces,
I.M.R.N. 1998, no. 19, 1011--1031.
\item{[B99]}
R. E. Borcherds,
The Gross-Kohnen-Zagier theorem in higher dimensions.
Preprint alg-geom/9710002. 
Duke Math. J. 97 (1999), no. 2, 219--233. 
\item{[Br]} J. H. Bruinier,
Borcherds products and Chern classes of Hirzebruch-Zagier divisors,
1999 preprint \hfill 
available from
http://www.rzuser.uni-heidelberg.de/\hbox{\~{}}t91/bruinier/bochern.dvi
\item{[B]}
V. O. Bugaenko, Arithmetic crystallographic groups generated by
reflections, and reflective hyperbolic lattices. Lie groups, their
discrete subgroups, and invariant theory, 33--55, Adv. Soviet Math.,
8, Amer. Math. Soc., Providence, RI, 1992.
\item{[C-S]}
J. H. Conway,  N. J. A. Sloane,
 Sphere packings, lattices and groups. Second
edition. Grundlehren der Mathematischen Wissenschaften
 290. Springer-Verlag, New York, 1993.
 ISBN: 0-387-97912-3
\item{[E94]} W. Ebeling,
Lattices and codes. A course partially based on lectures by F.
Hirzebruch. Advanced Lectures in Mathematics. Friedr. Vieweg \& Sohn,
Braunschweig, 1994.  ISBN: 3-528-06497-8
\item{[E]}
F. Esselmann \"Uber die maximale Dimension von Lorentz-Gittern mit
coendlicher Spiegelungsgruppe.  (German) [On the maximal dimension of
Lorentzian lattices with co-finite reflection group] J. Number Theory
61 (1996), no. 1, 103--144.
\item{[F]}
J. Fischer,   An approach to the Selberg trace formula via the
Selberg zeta-function. Lecture Notes in Mathematics,
1253. Springer-Verlag, Berlin-New York, 1987.  ISBN:
3-540-15208-3
\item{[F99]}
E. Freitag,
Some modular forms related to cubic surfaces,
1999 preprint, available from www.rzuser.uni-heidelberg.de/\hbox{\~{}}t91
\item{[Ko]}
N. Koblitz,
Introduction to elliptic curves and modular forms. Second edition.
Graduate Texts in Mathematics, 97. Springer-Verlag, New York, 1993.
ISBN: 0-387-97966-2
\item{[K-K]} J. H. Keum, S. Kondo, The automorphism  groups of Kummer surfaces
associated with the product of two elliptic curves, 1998 preprint.
\item{[K]} S. Kondo, The automorphism group of a generic Jacobian Kummer
surface. J. Algebraic Geometry 7 (1998) 589--609.
\item{[M]}
Y. Martin, 
Multiplicative $\eta$-quotients, 
Trans. A.M.S. 348 no. 12 (1996) 4825--4856. 
\item{[Mi]}
T. Miyake,  Modular forms.
Springer-Verlag, Berlin-New York, 1989. ISBN: 3-540-50268-8
\item{[N57]}
M. Newman, Construction and application of a class of modular functions. Proc.
London. Math. Soc. (3) 7 (1957), 334--350.
\item{[N59]}
M. Newman, Construction and application of a class of modular functions. II.
Proc. London Math. Soc. (3) 9 1959 373--387
\item{[Ni]}
P. Niemann,
Some generalized Kac-Moody algebras with known root multiplicities,
Cambridge PhD thesis, 1997.
\item{[N]}
V. V. Nikulin, Integer symmetric bilinear forms and some of their
geometric applications. (Russian) Izv. Akad. Nauk SSSR Ser. Mat. 43
(1979), no. 1, 111--177, 238. English translation in Mathematics of
the U.S.S.R., Izvestia, Vol 14 No. 1 1980, 103-167.
\item{[N97]}
V. V. Nikulin,
On the classification of hyperbolic root systems of the rank three.
Part I alg-geom/9711032,
Part II alg-geom/9712033,
Part III math.AG/9905150.
\item{[R]}
H. Rademacher,
Topics in analytic number theory.
Die Grundlehren der mathematischen Wissenschaften, Band 169.
Springer-Verlag, New York-Heidelberg, 1973.
\item{[Ri]}
K. A. Ribet,
Galois representations attached to eigenforms with Nebentypus.
Modular functions of one variable, V
(Proc. Second Internat. Conf., Univ. Bonn, Bonn, 1976), pp.
17--51. Lecture Notes in Math., Vol. 601, Springer, Berlin, 1977.
\item{[S-H]}
R. Scharlau, B.  Hemkemeier,
Classification of integral lattices with large class number.
Math. Comp. 67 (1998), no. 222, 737--749.
\item{[S-V]}
R. Scharlau, B. B.  Venkov,  The genus of the Barnes-Wall lattice.
Comment. Math. Helv. 69 (1994), no. 2, 322--333.
\item{[Sh]}
G. Shimura,
Introduction to the arithmetic theory of automorphic functions. Kan\^o
Memorial Lectures, No. 1. Publications of the Mathematical Society of Japan,
No. 11. Iwanami Shoten,
Publishers, Tokyo; Princeton University Press, Princeton, N.J., 1971.
\item{[S-S]}
J-P. Serre, H. M.  Stark,
Modular forms of weight $1/2$. Modular functions of one variable, VI
(Proc. Second Internat. Conf., Univ. Bonn, Bonn, 1976), pp. 27--67.
Lecture Notes in Math., Vol. 627, Springer, Berlin, 1977.
\item{[V75]}
\`E. B. Vinberg, 
Some arithmetical discrete groups in \hbox{Loba\v cevski\u\i}  spaces.
Discrete subgroups of Lie groups and applications to moduli
(Internat. Colloq., Bombay, 1973), pp.  323--348. Oxford Univ. Press,
Bombay, 1975.
\item{[V85]}
\`E. B. Vinberg,
The two most algebraic $K3$ surfaces. Math. Ann. 265 (1983), no.  1,
1--21
\item{[V-K]}
\`E. B. Vinberg,  I. M. Kaplinskaja,  The groups $O\sb{18,1}(Z)$ and
$O\sb{19,1}(Z)$. (Russian) Dokl. Akad. Nauk SSSR 238 (1978), no. 6,
1273--1275.
\item{[Y]}
K-I Yoshikawa,
1999 preprint (?).
\bye